%Basic AmsTeX article template

\input amstex

\documentstyle{amsppt}

\overfullrule=0pt

\define\blkbox{{\vrule height1.5mm width1.5mm depth1mm}}

\topmatter

\title\nofrills Knots with unique minimal genus Seifert surface and depth of knots\endtitle
\author  Mark Brittenham \endauthor

\leftheadtext\nofrills{Mark Brittenham}
\rightheadtext\nofrills{ Unique minimal genus Seifert surfaces and depth}

\affil   University of Nebraska \endaffil
\address   Department of Mathematics, Oldfather 810, Lincoln, NE 68588-0323  \endaddress
%\curraddr           \endcurraddr
\email   mbritten\@math.unl.edu \endemail
%\dedicatory       \enddedicatory
%\date                    \enddate 
%%\thanks   Research supported in part by NSF grant \# DMS$-$ \endthanks
%\translator         \endtranslator
%\subjclass            \endsubjclass

\abstract 
We describe a procedure for creating infinite families of hyperbolic knots,
each having
unique minimal genus Seifert surface which cannot be the sole compact 
leaf of a depth one foliation.
\endabstract

\endtopmatter

\document

\heading{\S 0 \\ Introduction}\endheading

Thurston [Th1] has shown that every compact leaf $F$
of a taut foliation $\Cal F$ of a 3-manifold $M$ has least
genus among all surfaces representing the homology
class of the surface, that is, it realizes the 
Thurston norm of that surface.
Conversely, Gabai [Ga1] has shown that every
Thurston norm minimizing surface for a non-trivial 
homology class is the
sole compact leaf of a {\it finite depth}
foliation of $M$. A foliation $\Cal F$ has finite
depth if the leaves of $\Cal F$ can be partitioned
into finitely many classes ${\Cal D}_0$, ${\Cal D}_1$,
$\ldots$ ,${\Cal D}_n$, where ${\Cal D}_0$ consists
of the compact leaves of $\Cal F$, and each leaf in
${\Cal D}_i$ limits only on leaves in 
${\Cal D}_0 \cup\cdots\cup{\Cal D}_{i-1}$. The smallest
$n$ for which this is true is called the depth of $\Cal F$.

These results have been successfully employed,
largely by Gabai, to compute 
the genera of many classes of knots, such 
as arborescent knots [Ga2] and the knots in 
the standard knot tables [Ga3], by 
constructing finite depth foliations of the knot
exteriors, with a 
candidate Seifert surface as sole compact leaf. 
For each of these constructions, the foliations built
have depth one. The smallest depth of a finite depth
foliation for a knot (with a Seifert surface as
sole compact leaf) is called the {\it depth} of the knot; 
all of these knots therefore have depth (at most) one.

Cantwell and Conlon [CC1] gave the first examples of
knots with arbitrarily high depth, by employing an iterated 
Whitehead doubling construction. The large number of 
non-parallel incompressible
tori in such a knot complement forces the depth of the knot to be
correspondingly high. They then asked whether or not every
hyperbolic knot, by contrast, must have depth one.

Kobayashi [Ko] showed that the answer to this question was
`No', by exhibiting a hyperbolic knot, having a unique 
minimal genus Seifert surface, which cannot be the 
leaf of a depth one foliation. It is, however, the leaf of a 
depth two foliation, and so Kobayashi's knot has depth two.

In this paper we show that Kobayashi's example is not alone;
we construct large numbers of 
hyperbolic knots with unique minimal genus Seifert
surface which cannot be the compact leaf of a depth one
foliation of the knot exterior. Our approach is to focus 
on knots with free Seifert surfaces of genus one, that is, Seifert surfaces
whose complements are handlebodies of genus two. By focusing
on surfaces of genus one, we can employ a necessary condition
for the surface to be a leaf of a depth one foliation, due 
to Cantwell and Conlon [CC2]; by using free Seifert surfaces,
we can check algebraically that their condition is not satisfied,
using Beiri-Neumann Strebel invariants [BNS] and a computational
tool due to Brown [Bro]. Finally, we can use some standard
cut and paste techniques, and the algebra of normal forms in a free group,
to show that many of the examples we construct have unique 
minimal genus Seifert surface.

\heading{\S 1 \\ The construction}\endheading

Our examples are based on the construction of 
knots with genus one free Seifert surfaces, given
in [Br1]. The basic idea is to start with a `base'
knot $K$ with genus one free Seifert surface
$F$, and then choose a simple loop $\gamma$ lying on the boundary of
$(S^3\setminus$int$N(F))$ = $X(F)$, which bounds a disk in $N(F)$; see
Figure 1. We choose curves $\gamma$ lying on the four-punctured sphere
shown there; such curves will always bound disks in $N(F)$. Because
$\gamma$ is unknotted in $N(F)$, hence in $S^3$,
$1/n$ Dehn surgery on $\gamma$ will yield $S^3$
back again. The knot $K$, and the Seifert surface $F$, 
are carried under the surgery to a new knot and Seifert surface $K_n$ and $F_n$
in $\gamma(1/n) = S^3$. By  [Br1, pp. 63-64], the surfaces $F_n$ are all free.

\input epsf.tex

\leavevmode

\epsfxsize=3in
\centerline{{\epsfbox{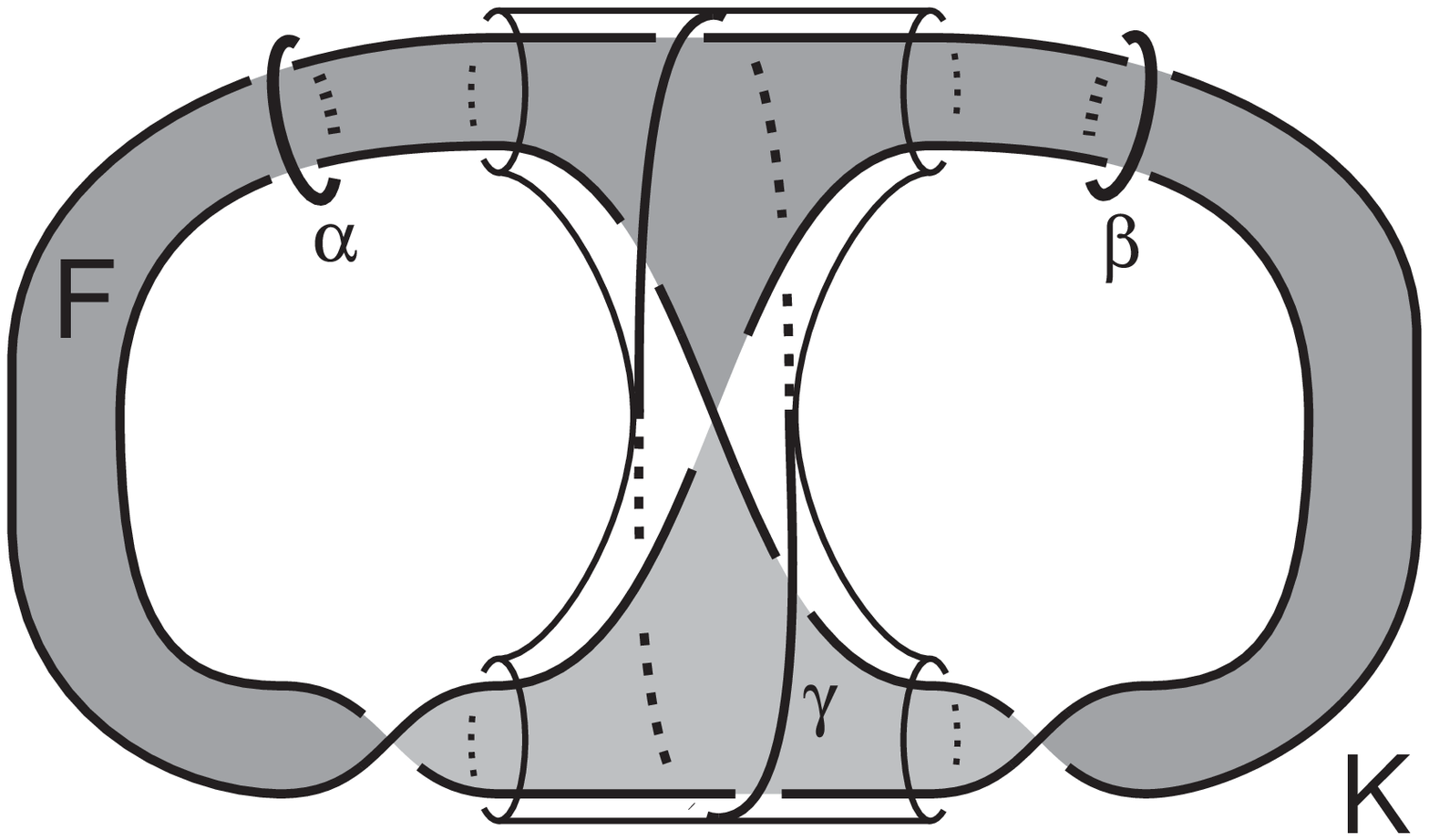}}}

\centerline{Figure 1}

\medskip

In anticipation of the next sections, to obtain the additional properties on 
$K_n$ and $F_n$ that we desire, we will also impose
three additional conditions on $\gamma$ (and $K$). 

First, $1/p$ and $1/q$ surgeries along the curves $\alpha$
and $\beta$ labelled in Figure 1 will simply add full twists to 
$K$, disjoint 
from the region of $\partial X(F)$ where we will be choosing $\gamma$,
and so will not much affect the construction. We can, in fact, think of 
ourselves as carrying out a \underbar{family} of constructions, one for each 
initial knot $K(1/p,1/q) = K_{p,q}$. For notational convenience, 
we will routinely suppress the
fact that $K_n$ really also depends on $p$ and $q$, that is, that the base knot
$K$ we are working with really depends on these two parameters. 
We need to keep these extra parameters in mind, however, because it will be large
numbers of twists (where here large essentially means $p,q\geq 2$)
which will, we shall see, insure that the resulting knots $K_n$ have
unique minimal genus Seifert surfaces. We assume the twists are
added so that they add to, rather than cancel out, the half-twists 
already present in $K$. Our first condition, therefore, is that our construction 
of the knots $K_n$ is actually built upon the base knots $K_{p,q}$.

Ultimately we will verify nearly all of the properties we need to establish 
about $K_n$ and $F_n$ by making computations in $\pi_1(X(F_n))$ = $F(a,b)$,
the free group on two letters $a$ and $b$.
This task will be much easier using a picture of $X(F)$ as a standard 
handlebody, with $K$, $\alpha$, $\beta$, and $\gamma$ drawn on its
boundary (where``$K$'' here really means a curve, isotopic to $K$, on $\partial X(F)$ 
which cobounds an annulus with $K$ in $N(F)$). By choosing the 
``obvious'' (vertical) compressing disks as 
the cores of our 1-handles in $X(F)$ to serve as the standard (horizontal) 
cores of the 1-handles for 
our handlebody, we obtain Figure 2 (see [Br1]). We have drawn
the curve $\gamma$ chosen in Figure 1, but a very large class of 
curves lying in the 4-punctured sphere complementary to $\alpha$ and 
$\beta$ will work for our purposes. In some sense, the best approach in 
fact is to find the necessary curve $\gamma$ in this standard picture, and then
determine which curve on $\partial X(F)\subseteq S^3$ it came from.

\leavevmode

\epsfxsize=5in
\centerline{{\epsfbox{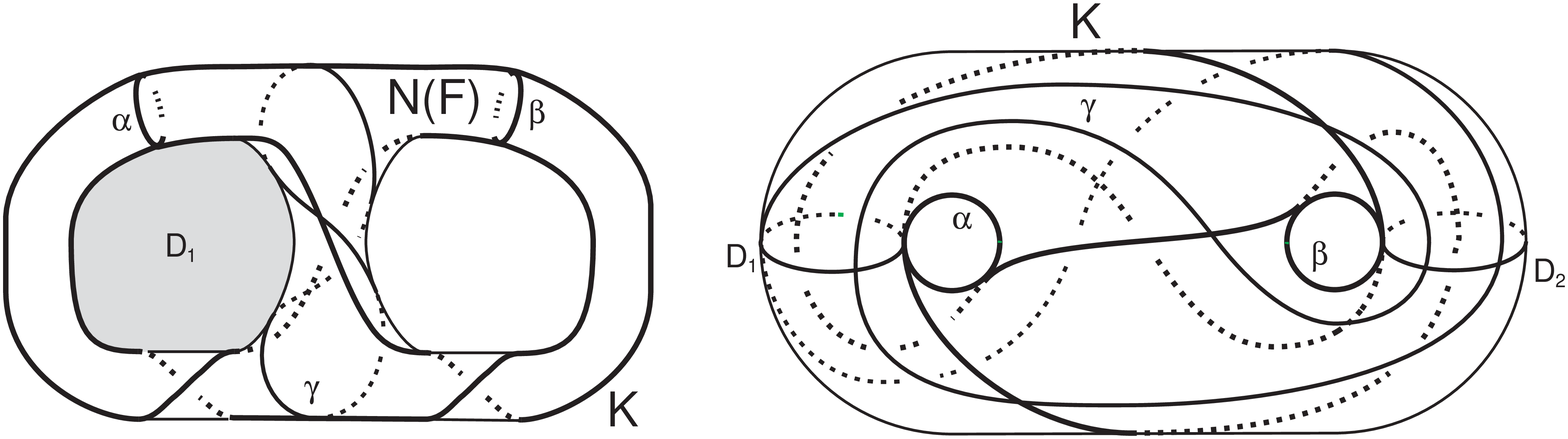}}}

\centerline{Figure 2}

\medskip

Throughout the rest of the paper, we will follow the standard practice 
of representing the inverses of the generators $a,b$ of the free
group $F(a,b)$ by $A$ and $B$, respectively, and the inverse of
a word $w$ in the letters $a,b$ by $\overline{w}$.

For our second condition, we choose a $\gamma$ which is homotopically essential in $X(F)$, 
but is null-homologous in 
$\partial X(F)$; that is,
$\gamma$ separates $\partial X(F)$ into two once-punctured tori. 
In terms of our standard
picture of $N(F)$ and $\gamma$, this is straightforward to
verify; we simply need to check that in the four-punctured
sphere in $\partial X(F)$ obtained by omitting
the two handles, in our figure, $\gamma$ does not separate 
the ends of either handle. This is sufficient, because $\gamma$ does 
separate the four-punctured sphere. This condition can also be readily 
established from the word in $F(a,b)$ representing $\gamma$; the
exponent sums of $a$ and $b$ in the word must both be zero.

When we straighten out the handlebody $X(F)$ to a 
standard handlebody, the knot $K$ and the loop $\gamma$
are carried to simple loops in $\partial X(F)$. We imagine
pushing $\gamma$ into $X(F)$ to a curve $\gamma^\prime$ on 
a surface parallel to $\partial X(F)$, to make it disjoint from $K$ again. 
The isotopies of $K$ (above) and $\gamma$
can be thought of as a single isotopy of $K\cup \alpha\cup\beta\cup\gamma$. 
Because $\gamma^\prime$ bounds a disk $D$ in $S^3$,
the effect of $1/n$ Dehn surgery on an object disjoint
from $\gamma^\prime$ is to cut the object open along $D$, 
give one side of the disk $n$ full twists, and then reglue.
As seen in Section 2 of [Br2], since the disk $D$ may be chosen to 
meet $X(F)$ in an annulus connecting $\gamma^\prime$ to $\gamma$, 
$X(F_n)$ and $K_n\subseteq \partial X(F_n)$
can be obtained from $X(F)$ and $K\subseteq \partial X(F)$ by 
simply setting $X(F_n)=X(F)$ and applying $n$ Dehn 
twists to $K$ along the curve $\gamma$, to obtain $K_n$ (Figure 3).

\leavevmode

\epsfxsize=3in
\centerline{{\epsfbox{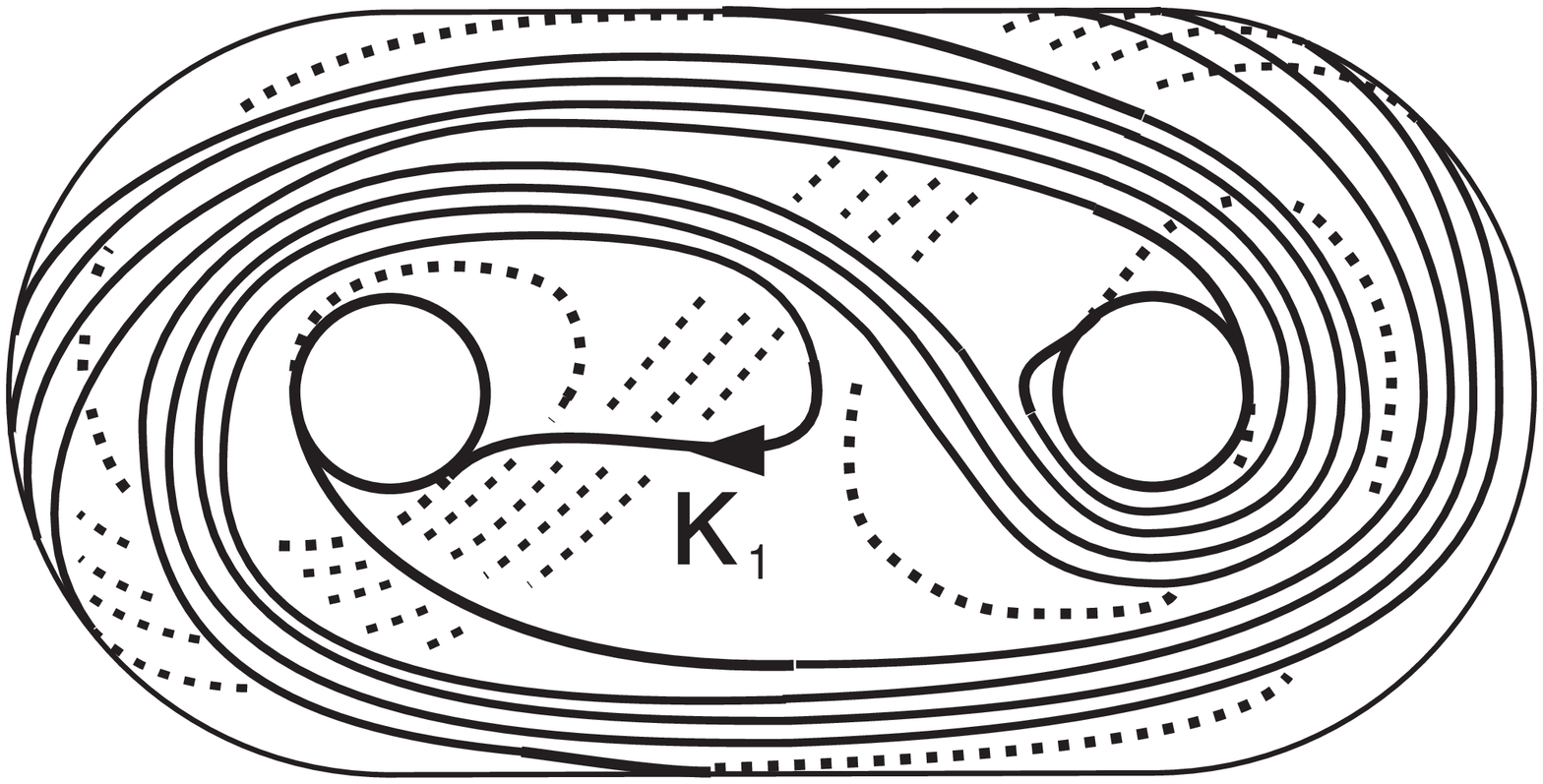}}}

\centerline{Figure 3}

\medskip

Because we start with a non-trivial knot $K$ - in our examples,
using $p,q\geq 2$, the knots $K$ we start with are
the 2-bridge knots [Schu] with continued fraction expansion
$[2p+1,-(2q+1)]$ - the surface $F$, being genus one, must be
incompressible, and so $\partial X(F)\setminus K$ is 
incompressible in $X(F)$. By applying a
criterion of Starr [Sta], we can also see this directly,
since there is a system of cutting disks for $X(F)$,
whose boundaries split $\partial X(F)$ into a pair of
thrice-punctured spheres $P_1,P_2$, so that $K$ meets each $P_i$ in 
essential arcs which join every pair of boundary curves.
In our case, we can use the standard set of disks 
$D_1,D_2,D_3$ where $D_1$ and $D_2$ are the disks chosen above, 
and $D_3$ is the horizontal disk in the middle of the 
handlebody in Figure 2b.

Our third condition on $\gamma$ is that it meet $K$ and the
cutting disks for $X(F)$ so that Dehn twisting $K$ along 
$\gamma^\prime$ yields a loop $K_1$ meeting the pair of 
thrice-punctured spheres in essential arcs. In particular,
this means that $\gamma$ itself cannot have a trivial 
arc of intersection with either of the three-punctured spheres $P_i$ (Figure 4).
Note that this condition is meant to apply to $\gamma$ 
{\it as drawn}, and not up to isotopy. This condition on $\gamma$
is not invariant under isotopy.

\leavevmode

\epsfxsize=2.4in
\centerline{{\epsfbox{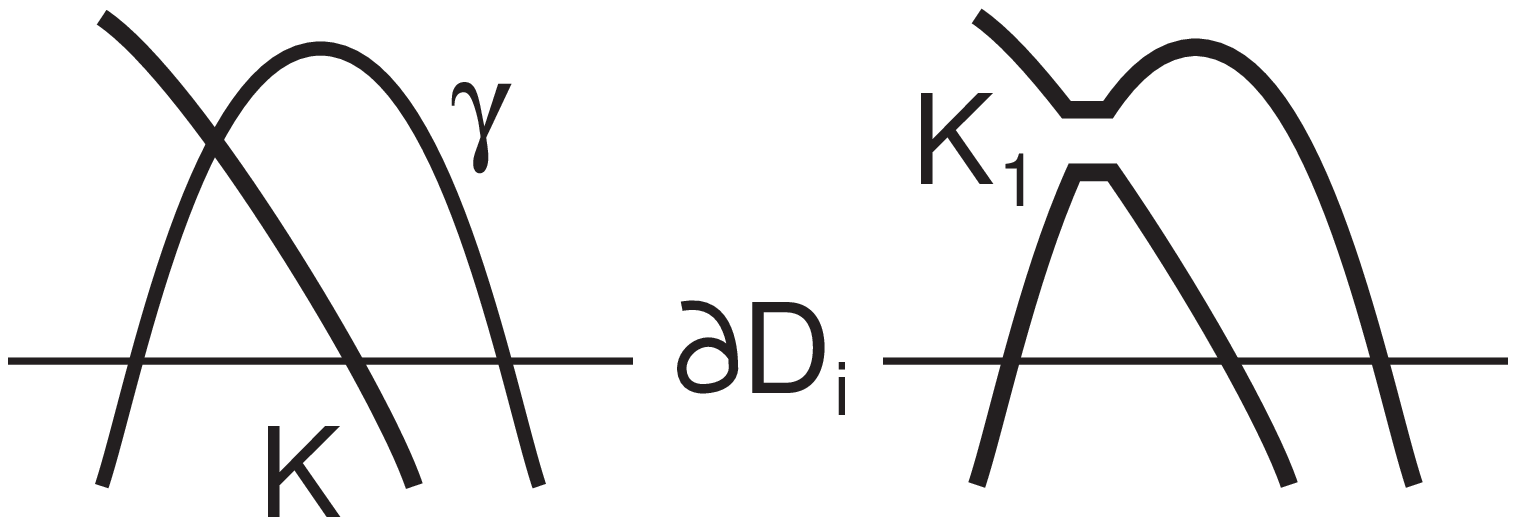}}}

\centerline{Figure 4}

\medskip

Since the arcs of intersection of $K_1$ with the $P_i$ basically 
include the arcs from 
$K$, $K_1$ meets each three-punctured sphere in all possible essential arc 
types, and so we immediately have that $\partial X(F_1)\setminus K_1$
is incompressible in $X(F_1)$, and so $K_1$ is non-trivial in $S^3$
and $F_1$ is a  least genus Seifert surface for $K_1$.

\leavevmode

\epsfxsize=5in
\centerline{{\epsfbox{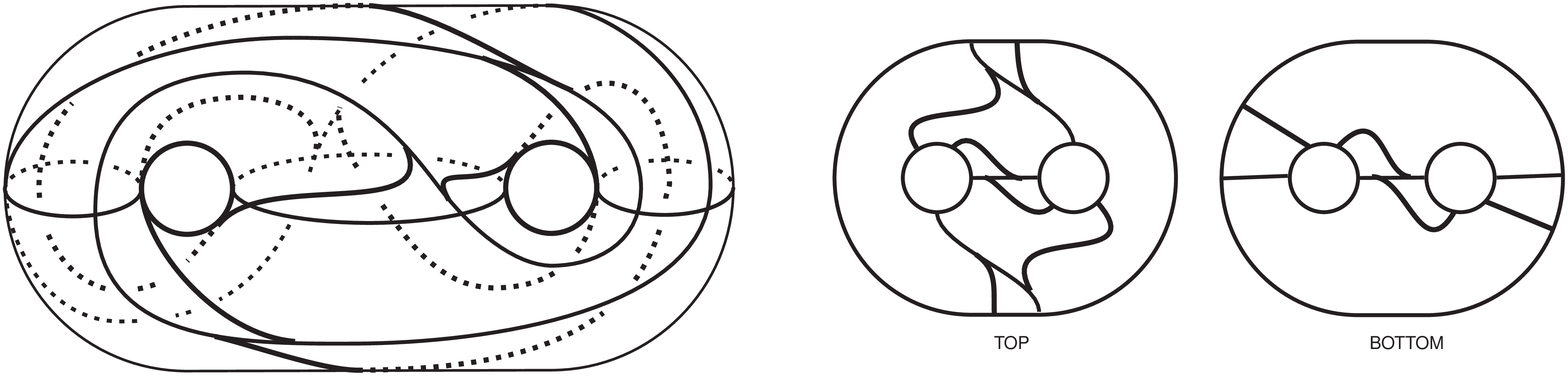}}}

\centerline{Figure 5}

\vfill\eject

This condition is best determined by building the train
track in $\partial X(F)$ which carries both $K$ and $K_1$; this amounts to 
turning each point of intersection of $K$ with $\gamma^\prime$ into a pair 
of switches (Figure 5).
We can then see directly that the intersection of the train track with
each three-punctured sphere does not carry a trivial arc, so long as it does
not carry an {\it outermost} one. This is essentially an Euler characteristic
argument. No horizontal boundary component of the train track is
boundary parallel, and no complementary region is a smooth disk or monogon.
Therefore, no complementary region has strictly positive (orbifold)
Euler characteristic, and so no union of regions can.

We provide several
examples of such $\gamma^\prime$ in Figure 6. Note that a Dehn twist in the opposite
direction yields a knot which does have trivial arcs, 
while twisting $n$ times in the direction in which one
twist works yields knots $K_n$ which also meet the punctured spheres in 
essential arcs, because the resulting curves are carried by the same train track as $K_1$.
$K_n$ is therefore non-trivial in $S^3$, and $F_n$ is least genus for $K_n$, for each $n\geq 1$.

\leavevmode

\epsfxsize=5in
\centerline{{\epsfbox{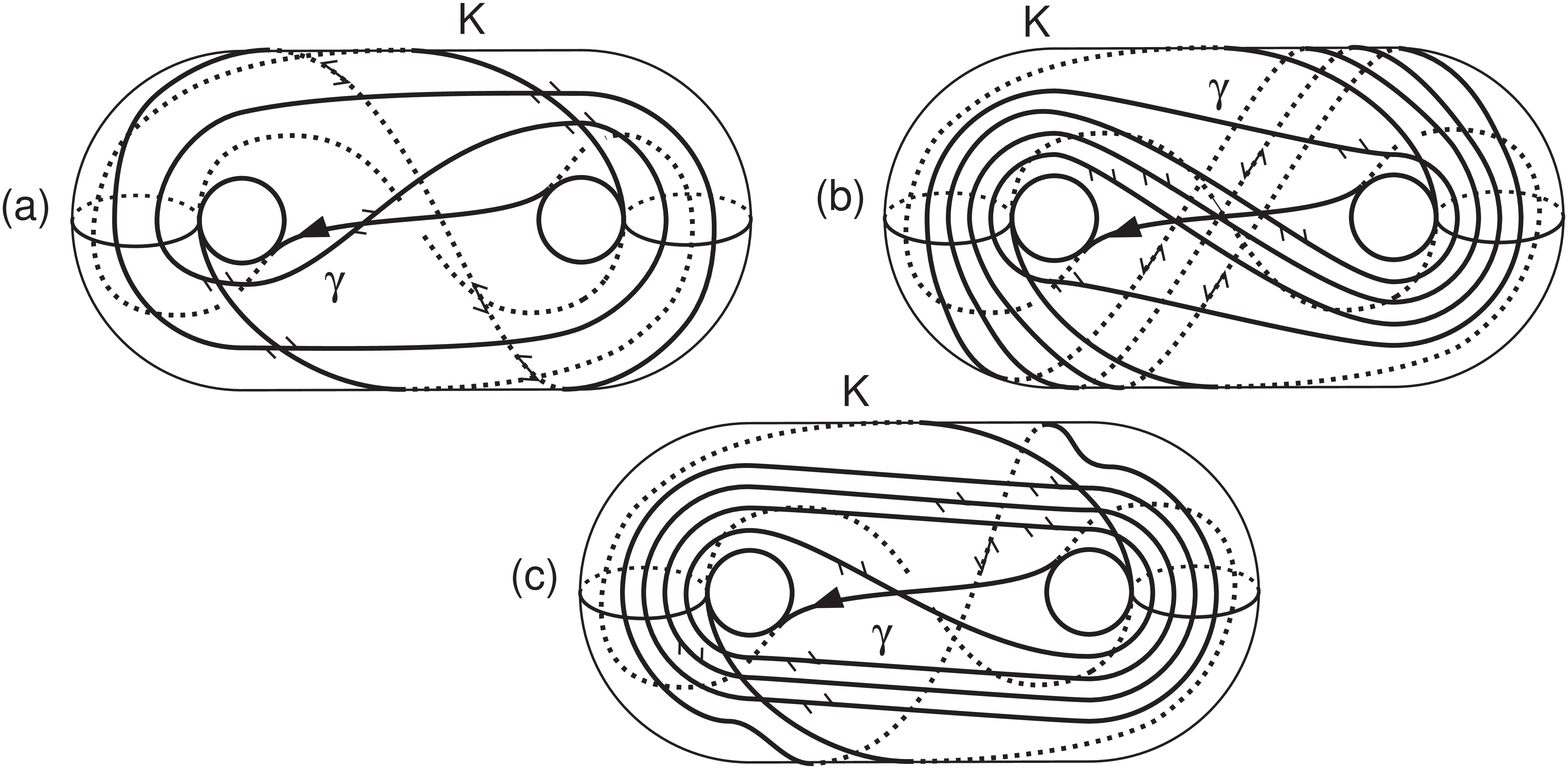}}}

\centerline{Figure 6}

\medskip

The fundamental group of $X(F)$ is free of rank two; it has a basis represented by loops, 
which we may take to lie in $\partial X(F)$, and which each intersect exactly once one of the two
standard compressing disks $D_1,D_2$ for $\partial X(F)$ labelled in Figure 2. 
The loop which intersects the
left disk, oriented to be travelling down as it passes through the disk, we will
denote $a$; the other, oriented similarly, will be denoted by $b$. An element $x$ of $\pi_1(X(F))$
can be written as a word in $a$ and $b$ and their inverses; this word can be read off from 
a loop representing $x$
by reading off the sequence of disks $D_1$ and $D_2$, and the directions, that $x$ passes through.

In the end, it will not matter much which specific words in $F(a,b)$ our knots $K_n$
spell out; our arguments turn out to be fairly general, relying only on the properties outlined
in Section 1. But with a little practice, it becomes a straightforward exercise to determine the word representing
the curve $K_n$, after $1/n$ surgery on $\gamma$, $1/p$ surgery on $\alpha$, and $1/q$ 
surgery on $\beta$. Note that an ordering for these surgeries does not need to be given;
because these curves bound disjoint disks in $S^3$, surgery on each is essentially
independent of the others. The basic idea is to first read off the word representing $\gamma$, 
as we traverse the curve. Then we read off the word representing $K_{p,q}$ (which is, in fact, $a^{p+1}b^{q+1}A^{p+1}B^{q+1}$), 
but each time we cross $\gamma$, we insert $n$ copies of a cyclic conjugate of $\gamma$ or 
$\overline{\gamma}$ into the word being read. Determining which cyclic conjugate is a
matter of bookkeeping, keeping track of where, as $\gamma$ is being read off, $K_{p,q}$ is crossing
$\gamma$. Note that the resulting word will already be in normal form in $F(a,b)$; there will be no occurance of a letter
followed by its inverse, because this would violate the condition that $K_n$ meet the $P_i$ in essential arcs.

For example, the words read off by the curves in Figures 3, 6(a), 6(b), and 6(c) are, in order

\smallskip

\noindent $A^{p+1} (baBA)^n b^{q+1} (bABa)^n a^{p+1}(BAba)^n B^{q+1} (BabA)^n$

\noindent $A^{p+1} (BabA)^n (baBA)^n (bABa)^n b^{q+1} (aBAb)^n a^{p+1} (bABa)^n (BAba)^n (BabA)^n B^{q+1}$

$(AbaB)^n$

\noindent $A^{p+1} (babaBABA)^n b^{q+1} (bABABaba)^n (baBABAba)^n (babABABa)^n a^{p+1}$

$ (BABAbaba)^n B^{q+1} (BababABA)^n (BAbabaBA)^n (BABababA)^n$ 

\noindent $A^{p+1} (bAbaBaBA)^n (bABaBabA)^n (baBaBAbA)^n b^{q+1} (bAbABaBa)^n a^{p+1}$

$(BaBAbAba)^n (BabAbABa)^n (BAbAbaBa)^n B^{q+1} (BaBabAbA)^n$

\smallskip

where we have read each word starting at the black triangle, reading in the direction that the triangle indicates.

Note that because the curves $\gamma\subseteq \partial X(F)$ that we choose are disjoint from $\alpha$ and $\beta$,
the word we read off in $F(a,b))$ for $\gamma$ never has a repeated letter. This is because for a repeated letter to occur ($b$, say),
$\gamma$ must pass through the disks $D_2$, $D_3$, and $D_2$, in order. This has the effect of "trapping" $\gamma$ against 
$\beta$ (Figure 7), which then requires $\gamma$ to pass through $D_2$ in the opposite direction before passing through $D_1$, 
giving an inessential intersection of $\gamma$ with $P_2$, contradicting our third assumption. The other possible repeated letters are similar.
Note also that the cyclic conjugates are never inserted within one of the powers, for the same reason;
the high powers arise within small neighborhoods of $\alpha$ and $\beta$, which $\gamma$ does not meet. 

\leavevmode

\epsfxsize=1.5in
\centerline{{\epsfbox{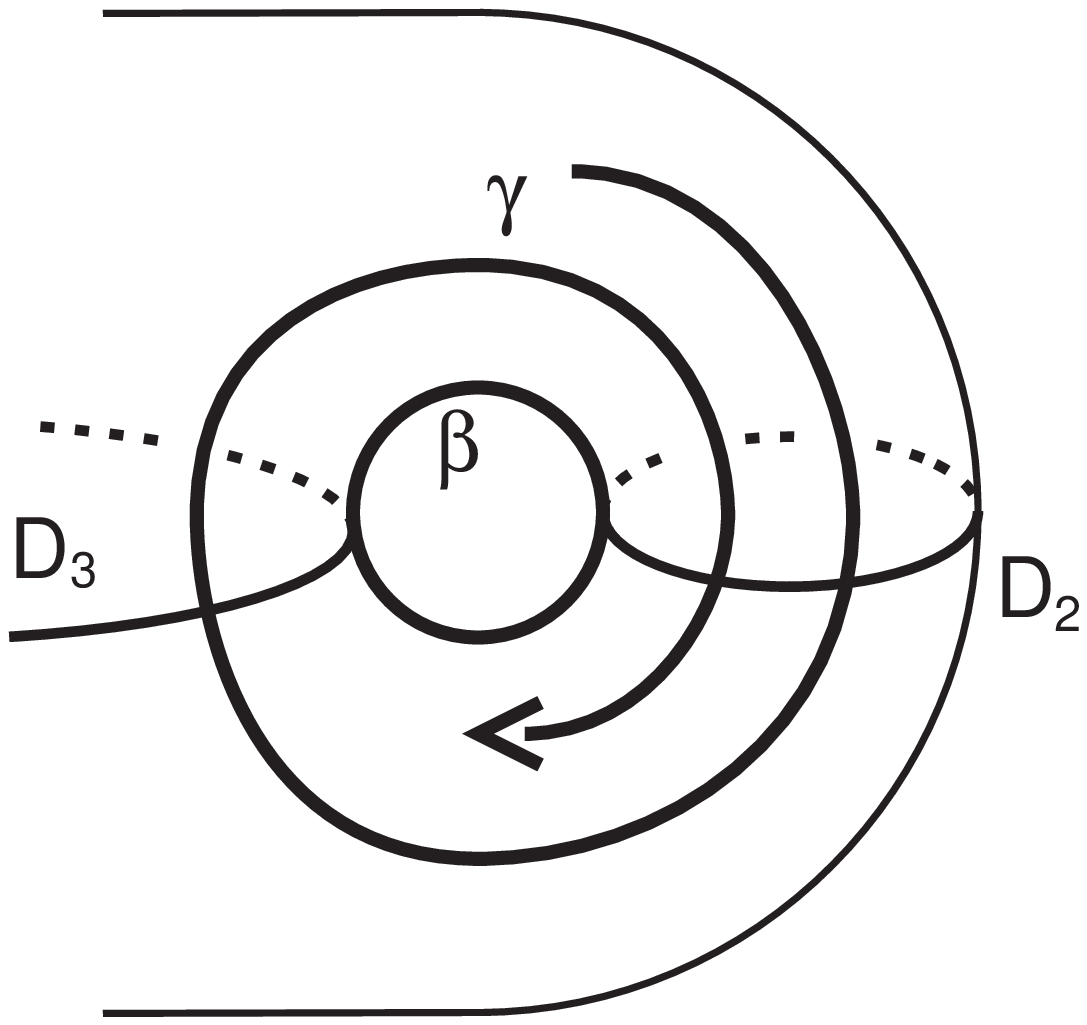}}}

\centerline{Figure 7}

\medskip

\heading{\S 2 \\ Depth greater than one}\endheading

We now show that for curves $\gamma$ satisfying the conditions above, the knot $K_n$ obtained 
by $1/n$ surgery along $\gamma$,
{\it for $n\geq 2$}, has free genus one Seifert surface $F_n$ which is not the sole compact
leaf of a depth one foliation. The key ingredient to showing this
is the following result, due to Cantwell and Conlon. Here $M$ denotes the manifold obtained
by gluing a 2-handle to $X(S)$ along the annulus $Q=X(S)\cap\partial(N(K))\subseteq \partial X(S)$

\proclaim{Theorem [CC2]}:
If $K$ is a genus-1 one knot and $S$ is a minimal genus Seifert surface for 
$K$, then every depth one foliation of $X(K)$, with $S$ as sole compact leaf, induces a fibration on
$M$, transverse to $\partial M$ and to a core of the attached 2-handle. Conversely, such a fibration
induces a depth-one foliation on $X(K)$ with $S$ as sole compact leaf.
\endproclaim

Therefore, to show that the Seifert surface $F_n$ for $K_n$ constructed as above is 
not the sole compact
leaf of a depth one foliation, it suffices to show that the manifold obtained
by gluing a 2 handle to $X(F_n)$ along $X(F_n)\cap\partial(N(K_n))\subseteq \partial X(F_n)$
is not a surface bundle over the circle. To show this, we will use the Bieri-Neumann-Strebel
invariant [BNS] of the fundamental group of $X(F_n)\cup_Q$(2-handle))= $M_n$. In the setting
of 3-manifold groups, the BNS invariant is essentially the same as the fibered faces of the
unit ball in the Thurston norm [Th1];  the BNS invariant is an open subset in the 
unit sphere of $H_1(M_n;{\Bbb R})$ = Hom$(\pi_1(M_n),{\Bbb R})$ whose points 
represent homomorphisms
of $\pi_1(M_n)$ to ${\Bbb Z}$ having finitely generated kernel. By the Stallings Fibration
Theorem [Stl], such a homomorphism is represented by a bundle map from $M_n$ to 
the circle (and conversely). Therefore, if the BNS invariant is empty, there are no 
homomorphisms from $\pi_1(M_n)$
to ${\Bbb Z}$ with finitely generated kernel, so $M_n$ is not a bundle over the circle.

$\pi_1(M_n)$ is, by the Seifert-Van Kampen Theorem, a one-relator group, since
$M_n$ is obtained from the genus-2 handlebody (with free fundamental group)
by adding a single 2-handle (giving the relator). The relator is the word, 
in the generators $a$ and $b$ of the free group, representing the knot $K_n$,
since $K_n$ is the core of the gluing annulus for the 2-handle. 

\smallskip

Brown [Bro] has produced an algorithm for computing the BNS invariant of a 1-relator group
$G$ = $<a,b$ $|$ $R>$.
His method requires that the relator $R$ have trivial abelianization; this will be the case for
our groups, since $K_n$ bounds a punctured torus in $\partial X(F)$, so is null-homologous.
$K_n$ is therefore trivial in the abelianization of  $\pi_1(X(F))$, since this group is canonically 
isomorphic to the first homology of $X(F)$. His algorithm 
also requires that the relator $R$ be reduced and cyclically
reduced (i.e, does not have an adjacent pair, and does not start and end 
with a pair, of letters $a,A$ or $b,B$). This is also the case in our setting, because
$K_n|(D_1\cup D_2\cup D_3)$ consists of essential arcs in $\partial X(F_n)|(D_1\cup D_2\cup D_3)$.
In our setting, the BNS invariant of $G$ will be a subset of the unit sphere in 
$H_1(M_n;{\Bbb R})$. $H_1(M_n;{\Bbb R})$ can be identified with 
Hom$(\pi_1(M_n),{\Bbb R})$ = ${\Bbb R}^2$,
since Hom$(\pi_1(M_n),{\Bbb R})$ = Hom$(F(a,b),{\Bbb R})$, because the relator $K_n$, 
being null-homologous, is automatically sent to $0$.

\leavevmode

\epsfxsize=4in
\hfill{{\epsfbox{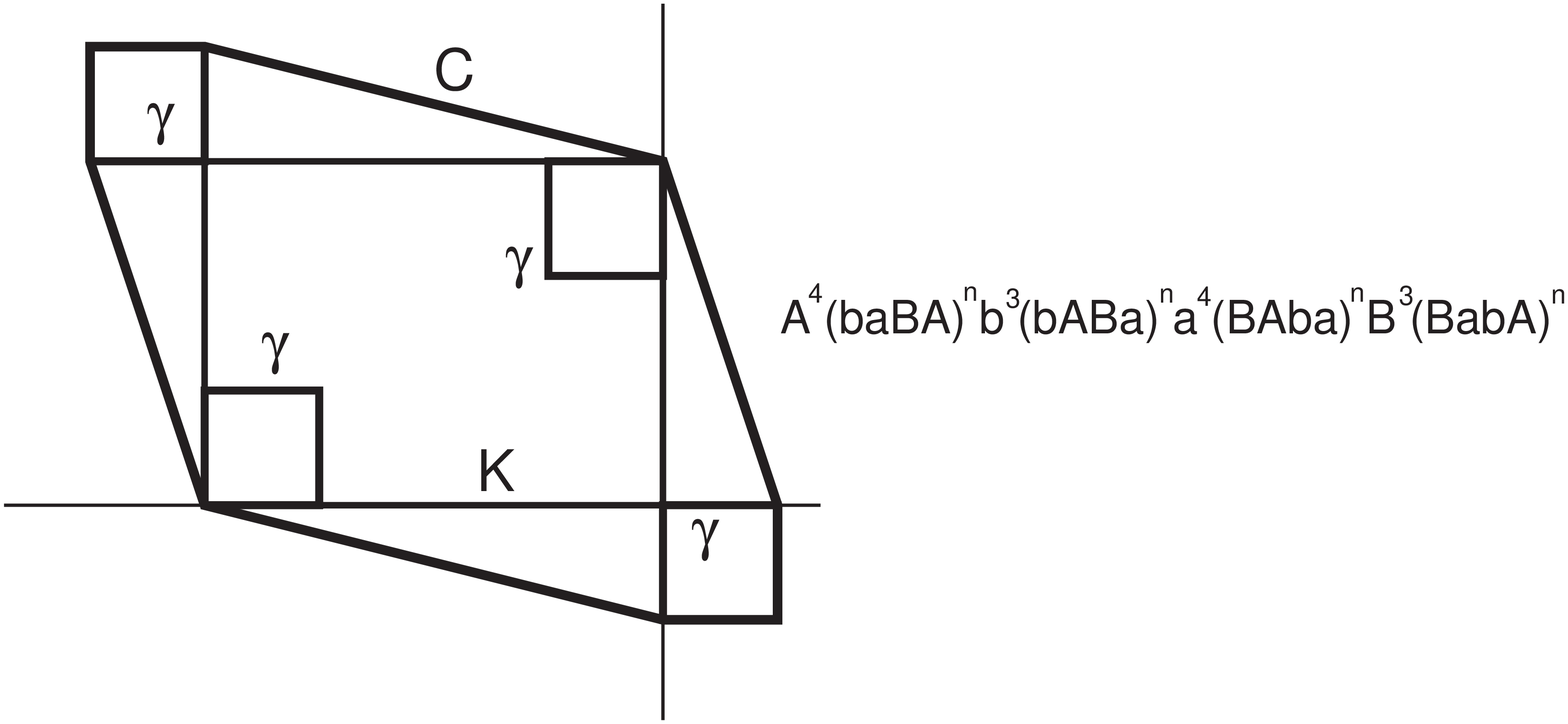}}}

\centerline{Figure 8}

\medskip

Brown's algorithm consists of tracing out the relator $R$ as a
path in ${\Bbb Z}^2\subseteq{\Bbb R}^2$, thinking of $R$ as lying in  ${\Bbb Z}^2$ instead of $F(a,b)$. 
However, we trace out $R$ in the order in which it is written; we do not try to 
simplify it in ${\Bbb Z}^2$ first; see Figure 8 for an example. Let $C$ denote the boundary of the convex hull
of the traced out path; it is a finite-sided convex polygon. A vertex $v$ of $C$ is called {\it simple}
if it is crossed by the path traced out by $R$ exactly once. $C$ will necessarily contain 
horizontal and vertical edges
(since the relator is cyclically reduced); such an edge is called {\it special} if both of its endpoints are
simple. Then [Bro, Theorem 4.4] the BNS invariant of $G$ consists of open arcs in the
unit circle in $H_1(M_n;{\Bbb R})$, one for each simple vertex or special edge in $C$,
{\it provided}, in the second case, that the line containing the special edge intersects $C$ only in 
that edge. His algorithm also describes how to compute these arcs.

For our purposes, however, we wish to establish that the BNS invariant is {\it empty}.
We therefore wish to find knots $K_n$ which, when thought of as words in the 
free group $F(a,b)$ and traced out as a path in ${\Bbb Z}^2\subseteq{\Bbb R}^2$,
have no simple vertices or special edges. That is, every vertex of $C$ is crossed at least twice by $K_n$.
For curves $\gamma$ and $K$ chosen as above, this is completely routine, provided we 
take $n\geq 2$, that is,
we spin $K$ at least twice around $\gamma$. The basic idea is that, since $\gamma$ is
null-homologous in $\partial X(F)$, and therefore in $X(F)$, when $K_n$ is traced out,
each occurrence of the cyclic conjugates of $\gamma$ and $\overline{\gamma}$ traces out a loop. Tracing each out
$n\geq 2$ times therefore insures that each vertex that these pieces of $K_n$ meet have
already been crossed at least twice by $K_n$. In fact, since with our choice of $K$
above $K_n$ is represented by a word of the form 
$A^{p+1}$(loops)$b^{q+1}$(loops)$a^{p+1}$(loops)$B^{q+1}$(loops), the path traced out
is essentially a large (for $p,q$ large) rectangle, with extra loops traced at the four corners.
Since every extra loop is traced out at least twice, the only possible simple vertices or
special edges that can be created are in the sides of the large rectangle, where any simple vertex will in 
fact be contain in a special edge. But the line containing such an edge will meet $C$ at the 
corners of the large rectangle, where the path has crossed at least twice. So neither criterion 
of Brown's algorithm can be met, and so all of the groups $<a,b$ $|$ $K_n>$ constructed in this way
have trivial BNS invariant.

Therefore, for every knot $K_n$ and free genus one Seifert surface $F_n$
constructed in this way, the 
manifold obtained by gluing a 2 handle to $X(F_n)$ along 
$X(F_n)\cap\partial(N(K_n))\subseteq \partial X(F_n)$
is not a surface bundle over the circle, provided that $n\geq 2$.
Therefore we have:

\proclaim{Proposition 1} For curves $\gamma$ satisfying the conditions of Section 1,  and for every $n\geq 2$,
the surface $F_n$ cannot be the sole compact leaf of a depth one foliation of $X(K_n)$.
\endproclaim

\heading{\S 3 \\ Unique genus one Seifert surface}\endheading

We now give a procedure for verifying that a Seifert surface like the ones
constructed above is the unique minimal genus spanning surface for the knot.
As in the previous section, the basis for the procedure is geometric, but
we will ultimately verify it algebraically, using the the word in the free group
$F(a,b)$ that the knot $K_n$ spells out.

Our starting point is a result of Scharlemann and Thompson, which says,
essentially, that to rule out a second non-isotopic minimal genus spanning
surface, we need only rule out the existence of a second such surface 
{\it disjoint} from the first:

%\vfill\eject

\proclaim{Theorem [ST]}: If $S$ and $T$ are minimal genus Seifert surfaces for the knot $K$,
then there is a sequence of minimal genus Seifert surfaces $S=S_0,S_1,\ldots,S_n=T$ such 
that, for each $i, 1\leq i\leq n, |S_i\cap S_{i-1}|=K$.
\endproclaim

So if one of the knots $K_n$ constructed as above has a second genus-1 Seifert surface $\Sigma$,
then it would have one which is disjoint from our surface $F_n$. $\Sigma$ is therefore
contained in the handlebody $H$ = $X(K_n)|F_n$ = $X(F_n)$, with boundary 
$\partial\Sigma = K_n \subseteq \partial H$. 
Because $\Sigma$ is not isotopic to $F_n$ in $S^3$, $\Sigma$ is not boundary parallel
in $H$. After isotopy, $\Sigma\cap D_3\subseteq D_3$ consists of circles and arcs,
and all circles of intersection can be removed by a standard innermost circle argument,
using the incompressibility of $\Sigma$. The intersection must still be non-empty,
because $\partial\Sigma = K_n$ meets $\partial D_3$. Choose an outermost arc $\eta$ 
of $\Sigma\cap D_3\subseteq D_3$, which together with an arc $\omega$ of $\partial D_3$
bounds a disk $\Delta\subseteq D_3$ whose interior is disjoint from $\Sigma$. $\eta$
cannot be $\partial$-parallel in $\Sigma$. For otherwise when we look at the intersection of 
the disk $E$ cut off by $\eta$ (with boundary $\eta\cup\eta_0$) with the cutting disks 
$D_1$ and $D_2$, an outermost 
(i.e, furthest from $\eta$) arc of intersection will cut an arc off of $\eta_0$, lying in $K_n$,
whose endpoints lie in $D_i$, $i=1$ or $2$. But this contradicts our third condition on $K_n$.
$\eta$ therefore is non-separating on $\Sigma$, and so when we $\partial$-compress
$\Sigma$ along $\Delta$, we obtain an annulus $A_1$ (Figure 9).

\leavevmode

\epsfxsize=4.5in
\centerline{{\epsfbox{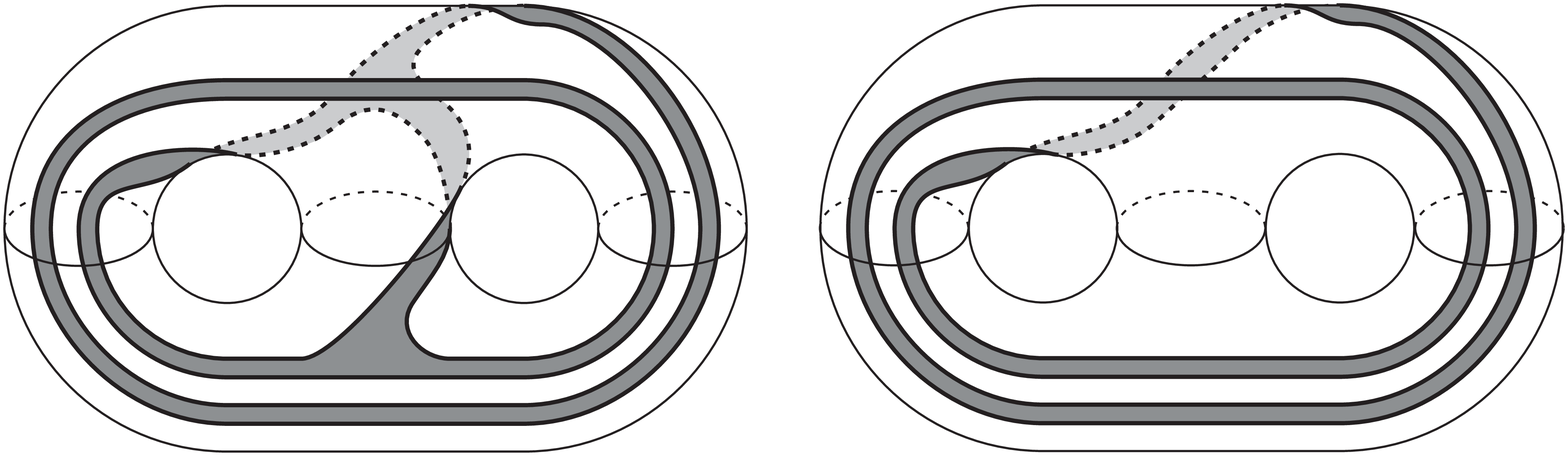}}}

\centerline{Figure 9}

\medskip

$\omega$ lies in one of the two punctured tori $T_1\cup T_2 = (\partial H))|K_n$, say $T_1$. 
$\omega$ does not separate $T_1$, because it is essential in $T_1$; the argument is identical 
to the one given above for $\Sigma$.
$T_1|\omega$ is therefore an annulus $A_2$ with 
(we may assume) $\partial A_2$ = $\partial A_1$ . $A_1\cup A_2= T\subseteq H$ is then a torus (it
cannot be a Klein bottle, because handlebodies do not contain non-orientable closed surfaces),
which must be compressible. Because $H$ is irreducible, $T$ either is contained in a 3-ball, so
the curves $c_1\cup c_2 = A_1\cap A_2$ are null-homotopic in $H$, or $T$ bounds a solid torus
$M_0$ in $H$. The first case is in fact impossible, because the curves $c_i$ are (non-separating, hence) essential in 
$T_1$, which is incompressible in $H$. In the second case, if $c_1$ (say) 
generates $\pi_1(M_0$), then $A_1$ is parallel to $A_2$
through $M_0$, that is, $A_2$ is $\partial$-parallel (Figure 10a). But this in turn implies that 
$\Sigma$ is $\partial$-parallel, since reversing the $\partial$-compression along $\Delta$, starting with a
$\partial$-parallel annulus, yields a $\partial$-parallel surface; $\Delta$ must lie outside of $M_0$,
since otherwise $\Sigma$ is compressible (Figure 10b).
$c_1$ must therefore represent a proper power of the generator of $\pi_1(M_0)$, and so represents a proper 
power in $\pi_1(H)$.

Therefore, for every outermost arc of $\Sigma\cap D_3\subseteq D_3$,
$\partial$-compression of $\Sigma$ along the outermost disk of $D_3$ that the arc cuts off produces an annulus
whose boundary components represent a proper power in the free group $F(a,b) = \pi_1(H)$. 

\leavevmode

\epsfxsize=3.5in
\centerline{{\epsfbox{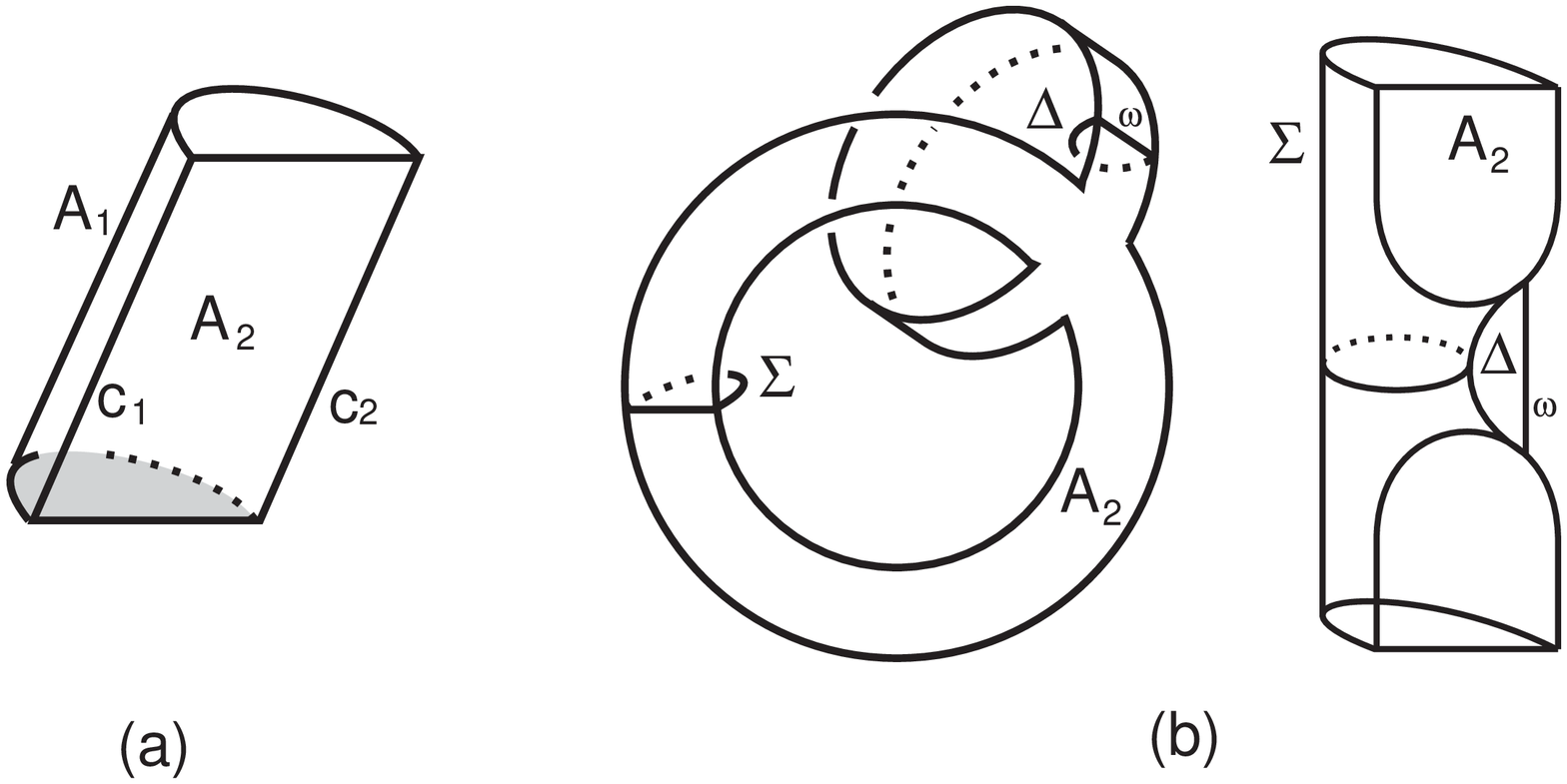}}}

\centerline{Figure 10}

\medskip

The main point, however, is that under this process of $\partial$-compression, curves representing 
only six words in $F(a,b)$ can be formed out
of any given curve $K_n$, and each can be quickly checked to determine if it is a proper power. The basic idea
is that, schematically, the punctured torus complements $T_1$ and $T_2$ of $K_n$ in $\partial H$ are built 
by gluing together pieces of 
$\partial H |(D_1\cup D_2 \cup D_3\cup K_n)$,
and $K_n$ cuts $\partial H |(D_1\cup D_2 \cup D_3)$ into four hexagons and a collection of rectangles, from the 
third condition on $\gamma$. 
The rectangles occur between arcs of $K_n$ parallel in $\partial H |(D_1\cup D_2 \cup D_3)=P_1\cup P_2$; two hexagons 
come from each of  $P_1,P_2$, lying between non-parallel arcs of $K_n$.
The $T_i$ therefore each look like a pair of hexagons with pairs of edges identified 
through a string of rectangles. From this point of view $T_i$ has a very natural spine consisting of two vertices, one for each 
hexagon, and three arcs, one for each string of rectangles (Figure 11). Each arc can be represented by a word 
$\lambda,\mu,\nu$ in $F(a,b)$ by reading off its intersections with the disks $D_1$ and $D_2$ dual to our generators $a,b$.
The word representing $K_n$, in terms of these three words, is then $\lambda\overline{\mu}\nu\overline{\lambda}\mu\overline{\nu}$.

\leavevmode

\epsfxsize=2.2in
\centerline{\hskip.7in{\epsfbox{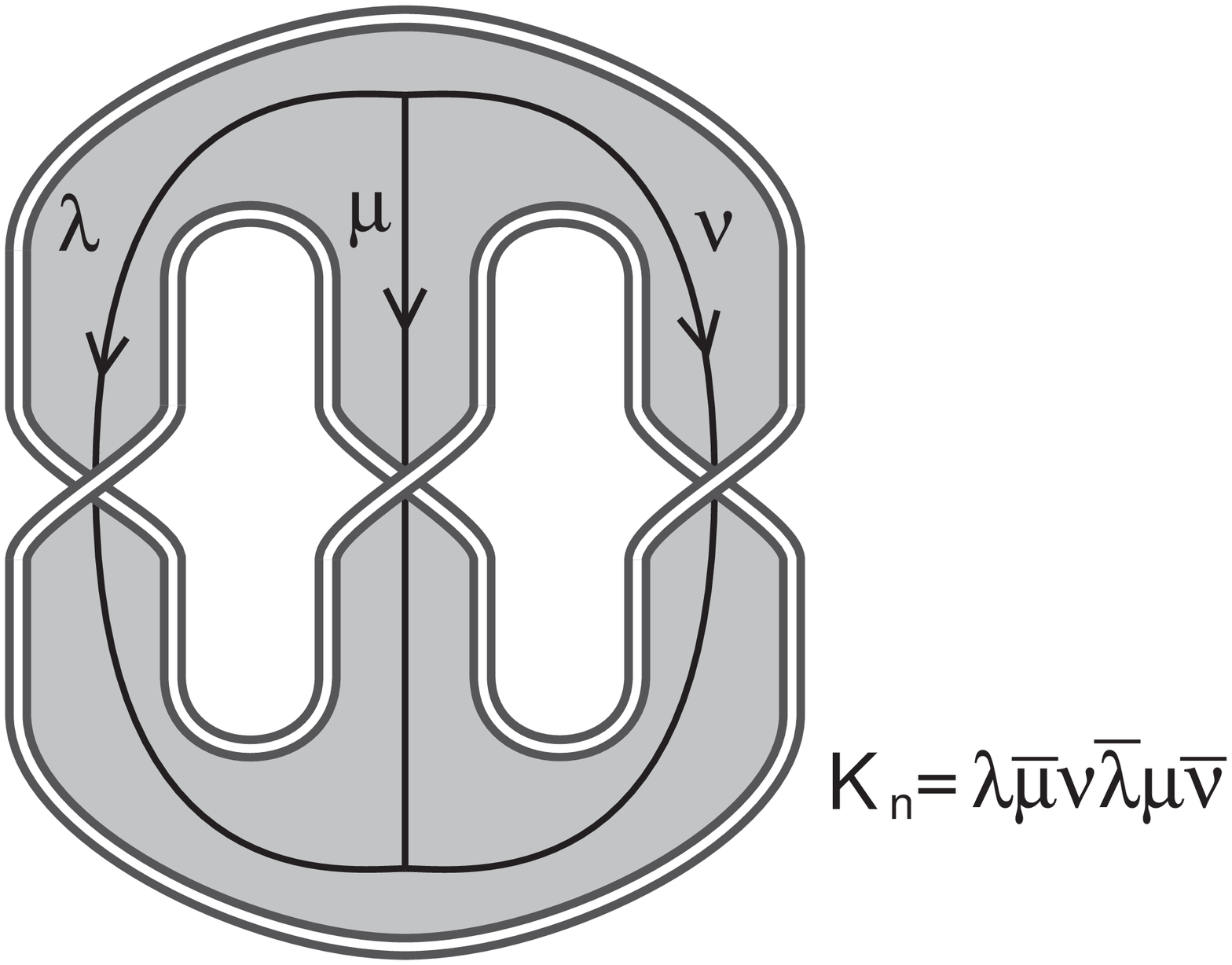}}}

\centerline{Figure 11}

\medskip

$K_n$ cuts $\partial D_3$ into a (potentially large) collection of arcs, and each of these arcs potentially cobounds an outermost
arc of $\Sigma\cap D_3$. We must therefore check that the annulus obtained by surgering $T_1$ or $T_2$, whichever one
contains the arc, along each arc does not have boundaries representing proper powers. But under
$\partial$-compression  each arc will cut across one of the three strings
of rectangles, resulting in an annulus made up, essentially, of just the other two strings. The core of this annulus will therefore be a
curve represented by one of the words $\lambda\overline{\mu}$, $\nu\overline{\lambda}$, or $\mu\overline{\nu}$, depending 
upon which string of rectangles was cut. We should note that these words will be (cyclically) reduced; if, for example,
$\lambda$ begins by passing through $D_1$ (so starts with $a$) and $\mu$ passes through $D_2$ (so starts with $b$), then 
$\nu$ must first pass through $D_3$ and so first passes through $D_1$ or $D_2$ in the opposite direction, and so (is trivial
or) starts with $A$ or $B$. Consequently, $\overline{\nu}$ {\it ends} with $a$ or $b$, so $\lambda\overline{\nu}$ and $\mu\overline{\nu}$
are cyclically reduced. They are, as subwords of $K_n$, already reduced. All other combinations of initial letters are similar.

In the end it will not be necessary to do so, but as a practical matter, finding the words $\lambda,\mu,\nu$ is fairly straightforward. 
It can be done from a picture of $K_n$, as in 
Figure 12 (for one arc in each of $T_1$ and $T_2$); once one of the words (for each punctured torus) is
determined, the others can easily be found, since the words found will typically occur only twice (once forward and once
backward) in $K_n$, and the remaining two pieces will have only one representation as $w_1\overline{w_2}$ and $\overline{w_1}w_2$.
In practice, in fact, only two representations of $K_n$ in the required form will typically exist at all.

\leavevmode

\epsfxsize=5in
\centerline{{\epsfbox{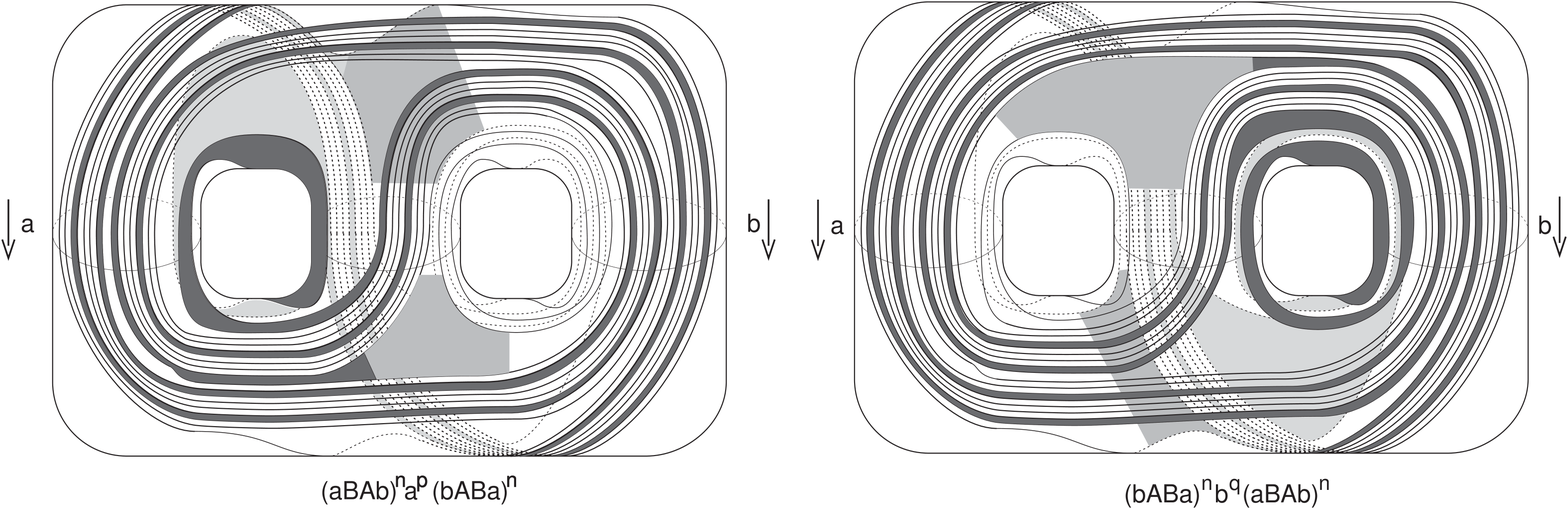}}}

\centerline{Figure 12}

\medskip

In the example of Figure 6(a), for example, we find that

\smallskip

\noindent $K_n$ = $[A^p][A(BabA)^n(baBA)^n][(bABa)^nb^{q+1}(aBAb)^n][a^p][a(bABa)^n(BAba)^n]$

\hskip1in $[(BabA)^nB^{q+1}(AbaB)^n]$ 

\noindent = $\lambda\overline{\mu}\nu\overline{\lambda}\mu\overline{\nu}$

\noindent and 

\noindent $K_n$ = $[(AbaB)^nA^{p+1}(BabA)^n][(baBA)^n(bABa)^nb][b^q][(aBAb)^na^{p+1}(bABa)^n]$

\hskip1in $[(BAba)^n(BabA)^nB][B^q]$ 

\noindent = $\lambda\overline{\mu}\nu\overline{\lambda}\mu\overline{\nu}$

\smallskip

\noindent representing the two different viewpoints of $K_n$ as $\partial T_i$.

As further examples, the reader can verify that in Figure 3 we have

$K_n$ = $[A^{p+1}][b][(aBAb)^nb^q(bABa)^n][a^{p+1}][B][(AbaB)^nB^q(BabA)^n]$

\hskip.2in so $\lambda = A^{p+1}$, $\mu = B$, and $\nu = (aBAb)^n b^q(bABa)^n$

$K_n$ = $[a][B^{q+1}][(BabA)^nA^p(AbaB)^n][A][b^{q+1}][(bABa)^na^p(aBAb)^n]$

\hskip.2in so $\lambda = a$, $\mu = b^{q+1}$, and $\nu = (BabA)^n A^p (AbaB)^n$

\medskip

Our goal now is to show that none of the resulting six  words (three for each of the $T_i$) 
$\lambda\overline{\mu}$, $\nu\overline{\lambda}$, or $\mu\overline{\nu}$ are proper 
powers in $F(a,b)$. To do this we use the fact
that the normal form of the $n$th power of a (cyclically reduced) word $w$ in a free group
 is the concatenation of $n$ copies of $w$. 
This is because the concatenation of $w$'s is, by our assumptions, already in normal form 
and cyclically reduced, and words in a free group have unique
normal form [MKS]. Conversely, if $v$ is cyclically reduced, and $v=w^n$, then the 
normal form for $w$ must be cyclically reduced (otherwise
$w^n$, when put into normal form, is \underbar{not}; the first letter of the first $w$ and 
the last letter of the last $w$ are inverses), and so 
the word $w^n$ is already in normal form, so $v$ is a concatenation of $w$'s. 
So we need only show that for each of the six words $w_j$ that the above process 
creates, $w_j$ is not a power of a subword of $w_j$.

This is where we will put the condition $p,q\geq 2$ to further use. The basic idea is that there 
is only one occurrence of  $a^{p+1}$ and $b^{q+1}$ in $K_n$,
and therefore at most one occurrence of a high power of $a,b,A$ or $B$, in each of the subwords 
$\lambda\overline{\mu}$, $\nu\overline{\lambda}$, or $\mu\overline{\nu}$ of $K_n$. Here ``high'' means ``$\geq 2$''.
By inspection, for $K_{p,q}$ itself, the subwords $\lambda, \mu, \nu$ are $a,A^p,b^{q+1}$ respectively, for one punctured torus, and
$a^{p+1},B^q,b$ for the other (Figure 13a).

The intersections of $\gamma$ with $K_{p,q}$ can be viewed (from the point of view of each of the $T_i$) as occurring in pairs, crossing the 
strings of rectangles for $T_i$ (since we can assume that $\gamma$  misses a pair of ``very thin'' hexagons in $T_i$ forming the 
complements of the strings of rectangles).
The effect on $\lambda, \mu, \nu$ of spinning $K_{p,q}$ around $\gamma$ is therefore to append
cyclic conjugates of $\gamma$ and $\overline{\gamma}$ before and after the values of $\lambda, \mu, \nu$ for 
$K_{p,q}$ (Figure 13b). These conjugates are not
inserted within $\lambda, \mu$, and $\nu$; this is because $\gamma$ is disjoint from $\alpha$ and $\beta$,
and the high powers in the strings of rectangles for $\partial X(F_{p,q})\setminus K_{p,q}$ can be thought of as originating in 
small neighborhoods of $\alpha,\beta$, when $K$ is spun about them to produce $K_{p,q}$. These regions of the strings 
of rectangles, corresponding to $\lambda, \mu, \nu$, are consequently also disjoint from $\gamma$.

\leavevmode

\epsfxsize=5in
\centerline{{\epsfbox{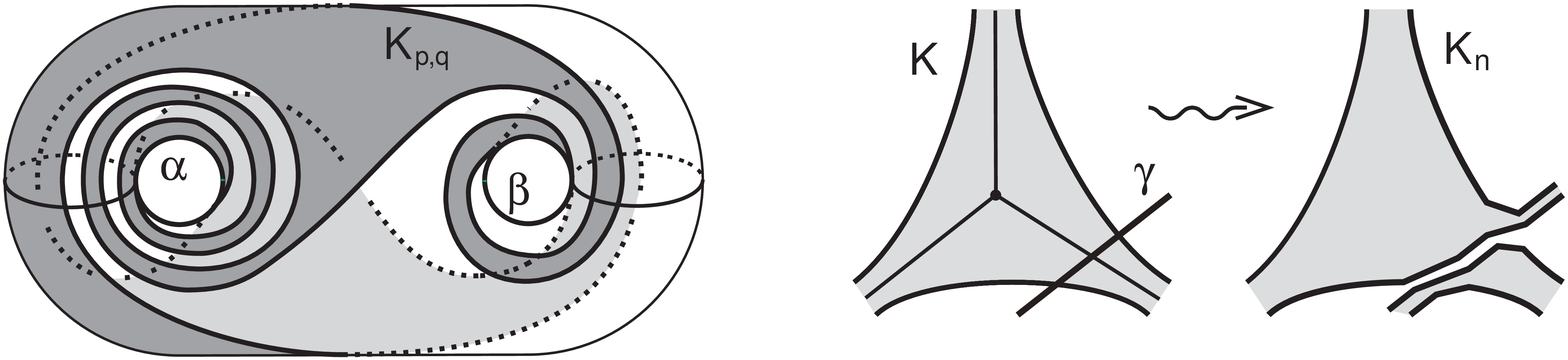}}}

\centerline{Figure 13}

\medskip

Each of the subwords $w$ = $\lambda\overline{\mu}$, $\nu\overline{\lambda}$ and $\mu\overline{\nu}$ for $F_{p,q}$
contains exactly one power $\geq 3$ of at least one of $a,b,A$ or $B$. For one of the two punctured tori, these words are 
$a^{p+1},b^{q+1}A$, and $aB^{q+1}$, respectively, and for the other are $a^{p+1}b^q, bA^{p+1}$, and $B^{q+1}$. 
The cyclic conjugates of $\gamma$ inserted into these words, used to obtain the corresponding words $w_j$, $j=1\ldots 6$, for $K_n$ 
contain no proper powers; therefore the only 
way that a word $w_j$ can be a proper power is if it contains one of these high powers and {\it no other letters}, because otherwise $w_j$ 
would be a power of a subword of $w_j$, and so would require more than one occurrence of this high power. 
This immediately rules out $\nu\overline{\lambda}$ and $\mu\overline{\nu}$ for the first punctured torus, and 
$\lambda\overline{\mu}$ and $\nu\overline{\lambda}$ for the second, as candidates for proper powers; their 
words already contain extra letters. 
The remaining two possibilities, $\lambda\overline{\mu}$ for the 
first torus and $\mu\overline{\nu}$ for the second, can be proper powers only if $\gamma$ does
not cross the pair of strings of rectangles making up one or the other of these two annuli. But the core of 
$(\lambda\overline{\mu})_1$, with word $a^{p+1}$,  intersects $\alpha$ once, and the core of 
$(\nu\overline{\lambda})_2$, with word $B^{q+1}$, intersects $\beta$ once (Figure 14).
So if $\gamma$ misses one of these cores $C$
then $\gamma$ lies in, and is null-homologous in (by one of our imposed conditions from section 1), 
the complement in $\partial X(F)$ of the corresponding pair of curves, which is a punctured torus. 
$\gamma$ is therefore isotopic to, and so freely homotopic to, 
the boundary of the punctured torus; the boundary is the only null-homologous essential simple curve 
in a punctured torus. $\gamma$ is then conjugate to the commutator of the two curves, which, since each is a power of the same generator $a$ or $b$, is trivial. So $\gamma$ must be null-homotopic,  contradicting 
another condtion on $\gamma$ from section 1. Consequently, we have:

\leavevmode

\epsfxsize=3in
\centerline{{\epsfbox{d_fig14.ai}}}

\centerline{Figure 14}

\medskip

\proclaim{Proposition 2} For every $p,q\geq 2$, and for every curve $\gamma$ 
satisfying the conditions of Section 1, the knots $K_n = K_{p,q,n}$ constructed above have unique minimal genus Seifert surface $F_n$; 
this surface is genus 1 and free.
\endproclaim

\heading{\S 4 \\ Hyperbolicity}\endheading

We now show that for each choice of $\gamma$ given in section 1, infinitely many of the resulting
knots $K_n$ are hyperbolic. To do this we demonstrate that for every $p,q\geq 2$ the two component link
$L$ = $K_{p,q}\cup\gamma$ has hyperbolic complement. Then by Thurston's Hyperbolic Dehn Surgery Theorem [Th2],
all but finitely many $1/n$ Dehn fillings along $\gamma$ yield hyperbolic knot complements. Since the volumes of the
complements of these knots $K_n$, for high values of $n$, are less than, but converge to, the volume of the complement of $K_{p,q}$
(whose volumes, in turn, converge to that of the three component link $K\cup\alpha\cup\beta$), infinitely many of the
knot complements are distinct. Therefore [GL] infinitely many of the knots are distinct.

Hyperbolicity can be verified by showing that the
topological hypotheses of Thurston's Geometrization Theorem [Th3] hold, that is, that $L$ is not a split, 
satellite, or torus link. The basic outline of the proof of this follows section 3 of [Br1]. We let $F_{p,q}$ denote the
genus 1 Seifert surface of $K_{p,q}$, obtained from $F$ by $1/p$ and $1/q$ surgery along $\alpha$ and $\beta$.

\medskip

\proclaim{Proposition 3} The link $L$ is not split.
\endproclaim

\medskip

{\bf Proof:} The loop $\gamma\subseteq X(F_{p,q})\subseteq X(K_{p,q})$ represents a
non-trivial word in $\pi_1(X(F_{p,q}))=F(a,b)$, by the second condition imposed on $\gamma$ above. Since $F_{p,q}$ is incompressible in 
$X(K_{p,q})$, $\pi_1(X(F_{p,q}))$ injects into $\pi_1(X(K_{p,q}))$,
and so $\gamma$ is non-trivial in $\pi_1(X(K_{p,q}))$. $\gamma$ is therefore not contained in a 3-ball in $X(K_{p,q})$. But then
for any 2-sphere $S^2$ in $X(L)\subseteq S^3$, whichever 3-ball complementary region (in $S^3$) of $S^2$ contains $K_{p,q}$ also contains 
$\gamma$. The other 3-ball region is therefore contained in $X(L)$. So $L$ is not split. \blkbox

\medskip

\proclaim{Proposition 4} The link $L$ is not a torus link.\endproclaim

{\bf Proof:} The knot $K_{p,q}$, since $p,q\geq 2$, is a non-torus alternating knot, and so by Menasco [Me] is hyperbolic. $L$
therefore has a non-torus knot component, and so is non-torus. \blkbox

\medskip

\proclaim{Lemma 5} There is no essential annulus $Q$ in $X(F_{p,q})\setminus$int $N(\gamma)$ = $M$ 
with one boundary component on $N(\gamma)$
and the other on $\partial X(F_{p,q})\setminus\partial N(K) = F_0\cup F_1$ .
\endproclaim

\medskip

{\bf Proof:} For any such annulus $Q$, $Q\cap\partial N(\gamma)$ = $\omega$ is a loop
representing some slope $u/v$ on the boundary torus (using the standard meridian/longitude coordinates). 
The other boundary component is on $F_i$, $i$ = 0 or 1. $v\neq 0$ since otherwise $Q$ capped off with a meridian
disk is a disk $D$ in $M$, with boundary in $F_1$, which is incompressible. $\partial D$
therefore bounds a disk in $F_1$, and so $D$ is $\partial$-parallel, by the irreducibility of $M$. In particular
$D$ separates $M$, but intersects the curve $\gamma$ once, a contradiction. 
Therefore $v\geq 1$ and $Q$ may then be capped of by an annulus in $N(\gamma)$ to (possibly a power of) the core of this solid torus, 
and so represents a homotopy between $\gamma^v$ and a simple loop on $F_1$. 

The manifolds $X(F_{p,q})$ 
are all handlebodies, as we have seen,  and we can in fact think of them as the \underbar{same} handlebody, 
with a different curve $K_{p,q}$ drawn on them,
obtained by spinning our original $K$, in our picture of $X(F)$, around $\alpha$ and $\beta$. 
From this point of view the curve $\gamma$ represents the same word in all of the free groups 
$\pi_1(X(F_{p,q})$, and via $Q$ is conjugate to a word in 
$i_*(\pi_1(F_0))$ or $i_*(\pi_1(F_1))$ in $\pi_1(X(F_{p,q}))$. But from Figure 13a we may 
compute that these two subgroups are freely generated by the words 
$\{a^{p+1},b^{q+1}A\}$ and $\{a^{p+1}B,b^{q+1}\}$, respectively. Consequently, 
if $p,q\geq 2$, then words in normal form in 
$i_*(\pi_1(F_0))$ always have the letter $b$ appearing as a proper power; no word in normal form 
can have a single $b$ or $B$ surrounded by a 
combination of $a$'s and/or $A$'s. Similarly, no word in normal form in 
$i_*(\pi_1(F_1))$ can have a single $a$ or $A$ surrounded by a combination of $b$'s and/or $B$'s. 
But for any curve $\gamma$ satisfying the conditions of Section 1, the word representing $\gamma$, by the comment at the end of Section 1, will
be cyclically reduced and will not have an 
occurrence of $a^2,A^2,b^2$, or $B^2$ in any cyclic conjugate. The same will therefore be true of $\gamma^v$, and so $\gamma^v$ cannot be conjugate 
to a word in either of these two subgroups.  \blkbox

\medskip

\proclaim{Proposition 6} $L$ is not a satellite link.\endproclaim

\medskip

{\bf Proof:} The argument parallels much of the argument of Proposition 3 of [Br1]; where the details are identical, we omit them. For ease of notation, 
we will write $K$ for $K_{p,q}$ and $F$ for $F_{p,q}$. In the discussion above, $\gamma$ was portrayed as lying on $\partial X(F)$; in the present discussion,
we think of it rather as lying slightly in the interior of $X(F)$.

\smallskip

Suppose that $T$ is an incompressible torus in $X(L)$ that is not $\partial$-parallel. Then since $F\subseteq X(L)$ is incompressible,
after disk swapping we may assume that $T\cap F\subseteq F$ consists of two families of parallel loops, one parallel to $\partial F$,
and one parallel to a non-separating loop in $F$. Repeatedly pushing $T$ across the annulus between the outermost loop of the first family and 
$\partial F$ yields a new essential torus (using Lemma 5 above in place of Lemma 1 of [Br1]), which we still call $T$;
after these isotopies, we may assume that the first family of loops is empty. All loops of $T\cap F$ are
then parallel to a single non-separating loop on $F$. 

Suppose that we still have $T\cap F\neq \emptyset$. For homological reasons ($T$ separates $X(K)$) 
there must be an even number of loops of $T\cap F$, cutting 
$F$ into an odd number of annuli $A_1,\ldots A_k$ and a single punctured annulus $P$. Because $X(K)$ is hyperbolic, 
$T$ must be compressible or $\partial$-parallel in $X(K)$.
It cannot be $\partial$-parallel, because then $P$ would lie in the product region between $T$ and $\partial X(K)$, so $P$ 
cannot $\pi_1$-inject into the product, implying that $F$ is compressible.
Since $T\subseteq X(K)\subseteq S^3$, $T$ bounds a solid torus $M_0$ in $S^3$ on at least one side. $T$ therefore either bounds the solid torus $M_0$ in $X(K)$,
if $K\cap M_0=\emptyset$, or $K\subseteq M_0$ and is disjoint from a meridian disk $D$ of $M_0$ (which is the only compressing disk $T$ has, unless 
$T$ bounds a solid torus on both sides). 

\leavevmode

\epsfxsize=2.3in
\centerline{{\epsfbox{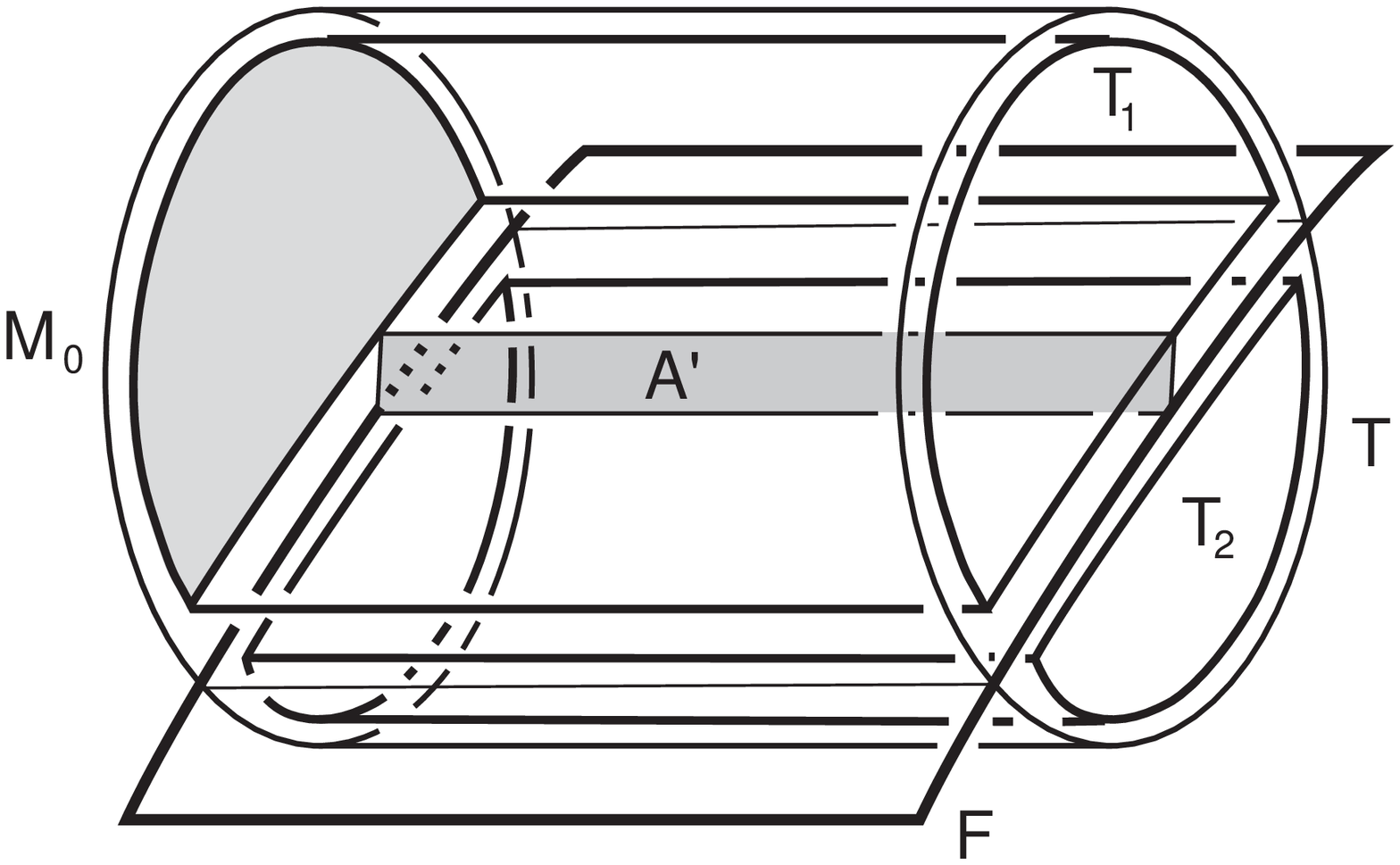}}}

\centerline{Figure 15}

\medskip

If $M_0\subseteq X(K)$, then at least one of the annuli $A_i$ lies in $M_0$; choose one that is outermost 
in $M_0$. Then $\partial N(T\cup A_i)$
consists of two tori $T_1$ and $T_2$, which lie in $X(L)$, since $\gamma$ is disjoint from $F$ 
(Figure 15). Because $A_i$ is outermost, one of these tori, say $T_1$, is 
disjoint from $F$. The argument of Proposition 3 of [Br1] applies without change to show that 
$T_1$ cannot be $\partial$-parallel in $X(L)$, and if 
$T_1$ is compressible in $X(L)$, then it bounds a solid torus in $X(L)$, and the solid torus may 
be used to isotope $T$ in $X(L)$ to reduce the number of curves 
of $F\cap T$. Therefore either $T_1$ is essential in $X(L)$ and disjoint from $F$, or we can 
reduce the number of circles of intersection of $T$ with $F$.

If $K\subseteq M_0$, then the argument of Proposition 3 of [Br1] applies to show that $T\cap F$ 
has more than two components, so one of the annuli
$A_i$ lies in $M_0$. Then we may repeat the argument of the previous paragraph to find either 
an essential torus $T_1$ disjoint from $F$ or reduce
the number of intersections of $T$ with $F$.

\leavevmode

\epsfxsize=3.2in
\centerline{{\epsfbox{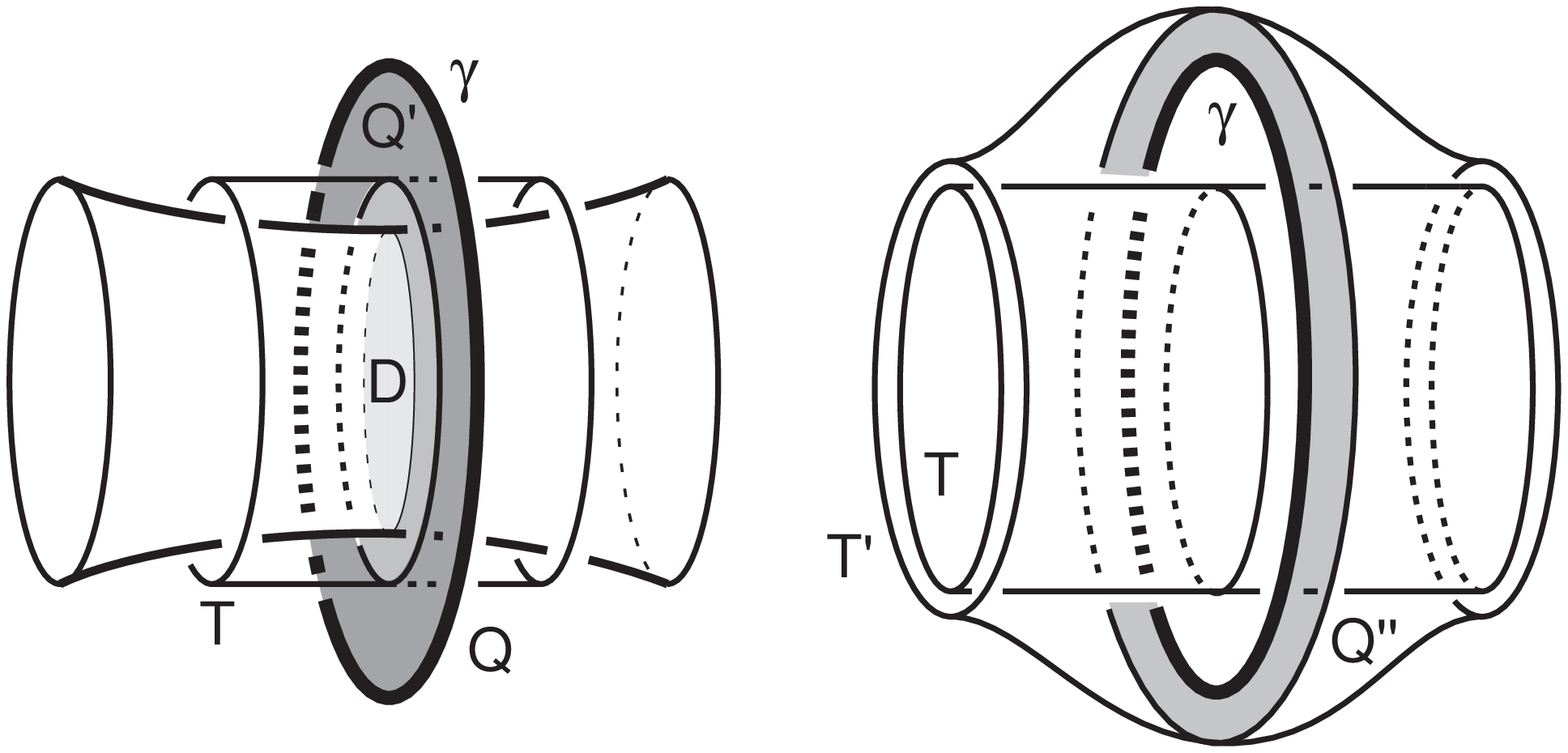}}}

\centerline{Figure 16}

\medskip

Therefore, by repeatedly applying these arguments, we eventually obtain an essential torus in $X(L)$ (which we will still call $T$) disjoint from $F$.
We then look at the intersection of $T$ with the annulus $Q$ lying between $\gamma$ and $\partial N(F)$ (Figure 16) where $Q$ is the intersection with $X(F)$
of a disk $D$ in $S^3$ bounding $\gamma$. (Recall that $\gamma$ bounds a disk in $X(F)$.) After isotopy in $X(L)$, by the incompressibility of $T$ we may assume that the intersection consists of loops essential in $Q$.
If $T\cap Q\neq\emptyset$, choose an outermost loop of intersection, cutting off an annulus $Q^\prime$ from $Q$ with one boundary component equal to 
$\gamma$. By isotoping $T$ across $Q^\prime$ (and therefore across $\gamma$), we obtain a new torus $T^\prime$ in $X(L)$. 

\smallskip

If $T^\prime$ is parallel in $X(L)$ to $\partial N(\gamma)$, then pushing $T^\prime$ back across $\gamma$ via the annulus $Q^{\prime\prime}$ to 
$T$ demonstrates that $T$ bounds a solid torus in $X(L)$, a contradiction. If $T^\prime$ is parallel in $X(L)$ to $\partial N(K)$, then since 
$T^\prime\cap F = \emptyset$, $F$ lies in the product region between $T^\prime$ and $\partial N(K)$, implying that $F$ is compressible, a 
contradiction. If $T^\prime$ is compressible in $X(L)$ via a compressing disk $D$, then as before either $T^\prime$ bounds a solid torus 
$M_1$ in $X(L)$ or it bounds a solid torus $M_1$ containing $K$ or $\gamma$, and the curve is disjoint from a meridian disk for $M_1$.

If $M_1\subseteq X(L)$, then $Q^{\prime\prime}\cap M_1 = \gamma_1$ must represent a generator of $\pi_1(M_1)$. For otherwise 
either $\gamma_1$ is a meridian curve, and then $Q^{\prime\prime}$ together with a meridian disk is a disk in $X(L)$ with boundary $\gamma$,
implying that $\gamma$ is null homotopic in $X(K)$, a contradiction, or $\gamma_1$ represents a proper power of the core of $M_1$. 
But then since $M_1\cup Q^{\prime\prime}\subseteq X(F)$, this implies that $\gamma$ represents a proper power in $\pi_1(X(F))=F(a,b)$,
and so, when written in normal form as a word in $a$ and $b$, it is a proper power of a subword. 
$\gamma$, however,  is parallel, via the annulus $Q$, to a simple curve $\eta$ in $\partial X(F) = G$, which, by a previous argument, cannot
be a proper power in $\pi_1(G)$. This in itself is not a contradiction, since there are simple curves in $G$ which are proper powers in 
$\pi_1(X(F))$. But $\eta$ is disjoint from the loops $\alpha$ and $\beta$ of Figure 1 above, and this implies that when written in 
the generators $a,b,c,d$ of $\pi_1(G) = <a,b,c,d : (AcaC)(dbDB)=1 >$ portrayed below (Figure 17), it can be chosen to be a word in only $a$ and $b$, since 
$G\setminus (D^2\cup\alpha\cup\beta)$ deformation retracts onto the subgraph generated by $a$ and $b$ of the spine of $G\setminus D^2$. The 
inclusion-induced map from $\pi_1(G)$ to $F(a,b)$ is the obvious one which sends $c$ and $d$ to 1; therefore, $\eta$, when written as a word 
in $a$ and $b$ in $\pi_1(G)$, is represented by the {\it same} word in $F(a,b)$. Since this word, in $F(a,b)$, is hypothesized to be a proper power,
$\eta$ therefore must represent a proper power in $\pi_1(G)$, a contradiction. 
But then when we push $T^\prime$ back to $T$ via $Q^{\prime\prime}$, $\gamma$ becomes the core of the solid torus $M_1$,
implying that $T$ is parallel to $\partial N(\gamma)$, a contradiction.

\leavevmode

\epsfxsize=2.4in
\centerline{{\epsfbox{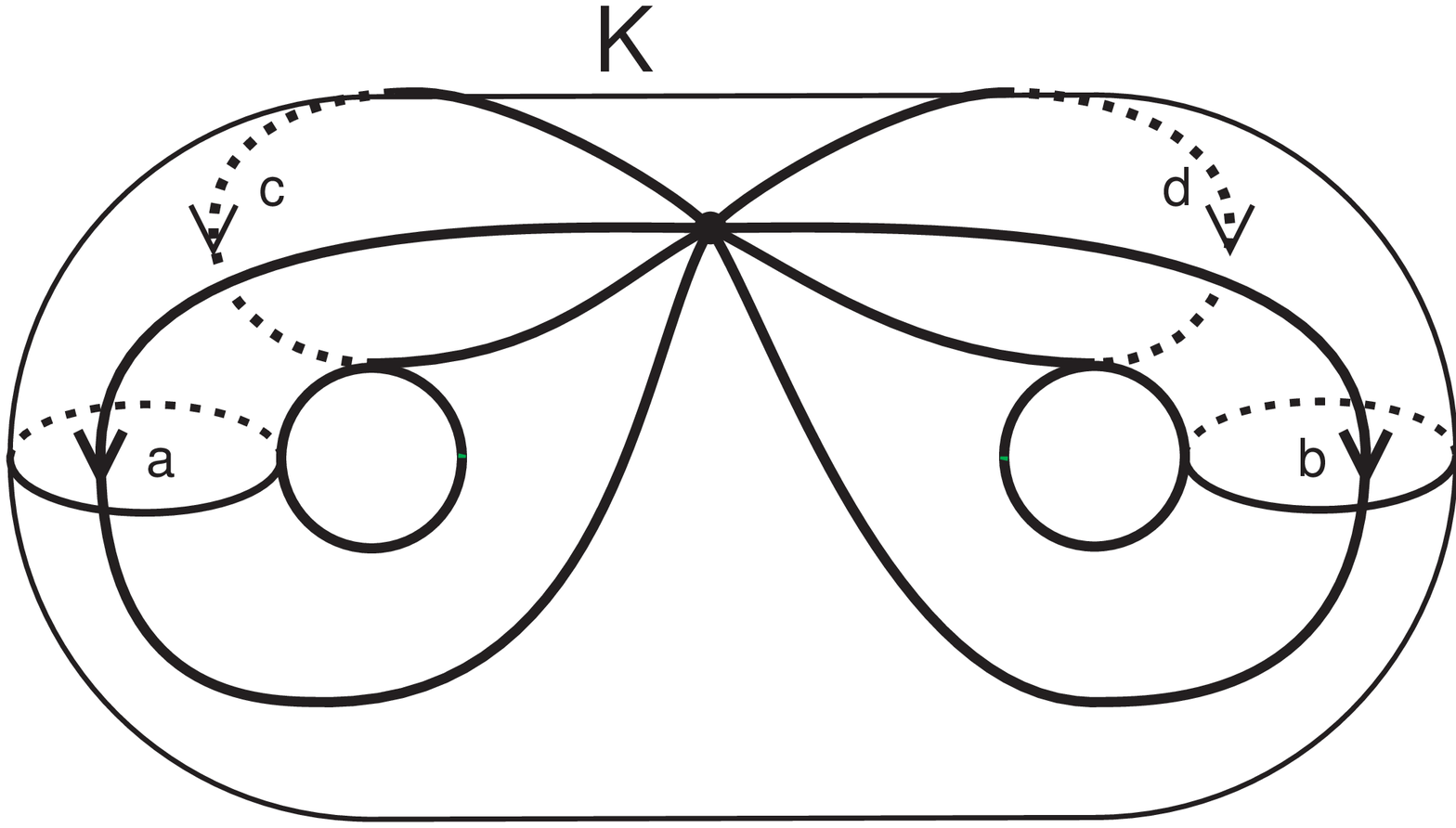}}}

\centerline{Figure 17}

\medskip

If  $T^\prime$  bounds a solid torus $M_1$ containing $K$ or $\gamma$, with this curve disjoint from a meridian disk $D$ for $M_1$,
then the argument of Proposition 3 of [Br1] applies to show that $D$ can be made disjoint from $Q^{\prime\prime}$. Then when we push $T^\prime$
back to $T$ via $Q^{\prime\prime}$, the compressing disk $D$ persists, so $T$ is compressible, a contradiction. 

Therefore, $T^\prime$ must be an essential torus in $X(L)$. Continuing this process of pushing the torus across $\gamma$, we eventually find a
torus $T$ essential in $X(L)$ which is disjoint from both $F$ and $Q$. But such a torus then lies in $X(F\cup Q)$, which is a handlebody, and so 
there is a compressing disk $D^\prime$ for $T$ in $X(F\cup Q)\subseteq X(L)$. $T$ is therefore compressible in $X(L)$, a contradiction.
So there are no essential tori in $X(L)$. \blkbox

\medskip

Consequently, the links $K_{p,q}\cup\gamma$, where $p,q\geq 2$ and $\gamma$ satisfies the conditions set forth in the previous sections, is a hyperbolic link.
Therefore,  for infinitely many $n$, the knots $K_n$ are hyperbolic.
Taking this and the results of the previous two sections together, we have:

\medskip

\proclaim{Theorem 7} There are infinitely many hyperbolic knots $K$ in the 
3-sphere whose exteriors do not admit a depth-one foliation.
\endproclaim

\heading{\S 5 \\ Concluding remarks}\endheading

In this paper we have found families of hyperbolic knots with depth greater than one. ``Found'' here must,
however be interpreted in the weak sense; by appealing to the Hyperbolic Dehn Surgery Theorem, we can 
conclude that our knots $K_n$ are hyperbolic for most values of $n$, but we can conclude this for
no {\it specific} value of $n$. Experimental evidence, by way of snappea [We] and Snap [Go], would suggest
that $K_n$ is in fact hyperbolic for ($p,q\geq 2$ and) all $n\geq 1$; sufficient diligence might allow for a 
unified proof of this. We do not undertake such an effort here.

Our result falls short of the standards set by Kobayashi's result in another way, in that we do not 
determine the actual depth of the knots $K_n$ that we have built. We only establish that it is at least 2, for each. Kobayashi
determined that the depth of his knot, which we will call $K_{ob}$, is exactly 2, by constructing a sutured manifold decomposition [Ga1]
of $X(K_{ob})$ by hand. We note that $K_{ob}$ is distinct from all of the knots we have built here, since the unique genus-1
Seifert surface for $K_{ob}$ does not have handlebody complement.

It is tempting to conjecture that all of the knots that we have built have depth 2,
since they are built from (2-bridge, hence) depth 1 knots by spinning around a loop $\gamma$ which bounds a surface $\Sigma$ in $X(F)$.
$\Sigma$ is also a natural candidate for the surface to next decompose the sutured manifold $(X(F_n),\partial X(F_n)\cap N(K_n))$ along.
The reader can verify, however, that this will not yield a taut sutured manifold; the surface $\Sigma$ can be chosen to be one
of the once-punctured torus components of $\partial X(F_n)\setminus\gamma$, and any choice of orientations yields trivial sutures
under the decomposition. Unfortunately, the technology does not yet exist to show that a hyperbolic knot has depth
greater than 2. To show that $X(K_n)$ does not admit a depth 2 foliation, we essentially need to show that
there is no decomposing surface $\Sigma$ for $X(F_n)$ such that decomposing $(X(F_n),K_n)$ along 
$\Sigma$, to $(M^\prime,\gamma^\prime)$,
yields a sutured manifold which admits a depth 1 foliation.
In our context, where $X(F_n)$ is a taut sutured handlebody, it may be possible to carry out such an analysis,
since every further decomposing surface $\Sigma$, being least genus and therefore incompressible,
will split $X(F_n)$ into another taut sutured handlebody; $M^\prime$
injects, on the level of fundamental group,
into the free group $\pi_1(X(F_n))$. We will therefore always remain within the realm of sutured handlebodies.

For our purposes, we can therefore focus on developing (more general) conditions to decide that a sutured handlebody
does not admit a depth 1 (or higher depth) foliation. For sutured handlebodies $(H,A)$ of any genus, depth 1
implies [Ga3] that int($H\cup_A(D^2\times I$) is a bundle over the circle, but when the genus of $H$ is greater than 2,
the fiber, necessarily, has infinite genus. When the genus of $H$ is 2, the fiber has finite genus, which is
why in the Cantwell-Conlon condition we may replace int($H\cup_A(D^2\times I)$ with $H\cup_A(D^2\times I)$ .
We can, at present, offer no insight into how to overcome the
problem of needing to test every possible decomposing surface.

Our techniques for establishing that $K_n$ had unique minimal genus Seifert surface
also used in a strong way the fact that the Seifert surfaces $F_n$ had genus 1. It 
would be interesting to develop more general conditions along these lines for 
generating knots with free genus $g$ Seifert surfaces which are unique for the knot
among minimal genus Seifert surfaces. In principle it seems likely that a similar approach,
using the group theoretic properties of the word representing $K_n$ in a free group,
should be successful. But the needed conditions will likely be more complicated.

Recent work of Scharlemann [Scha] implies that all of the
knots $K_n$ that we have constructed have tunnel number 2. A free genus one knot can have tunnel number 
at most 2; any two simple arcs in the surface $F_n$, which split $F_n$ to a 2-disk $D$, form a system of tunnels for $K_n$,
since the exterior of $K_n\cup$(tunnels) is $X(F_n)$, with the two disks in $F_+$ and $F_-$ corresponding to $D$ glued together.
But Scharlemann has shown that a genus-1 tunnel-1 knot must either be a satellite knot or a  2-bridge knot, which none of our 
examples are; 2-bridge knots are alternating, and so have depth 1 [Ga3].

\smallskip

In the end, it was surprising (to the author) how many of the original collection of knots $K_n$ had unique minimal 
genus Seifert surface. And so we end with the somewhat ill-posed

\proclaim{Question} How common is it for a knot to have unique minimal genus Seifert surface?
\endproclaim

\enddocument